\documentclass[12pt,draftcls,onecolumn]{IEEEtran}
\usepackage{cite}

\ifCLASSINFOpdf
  \usepackage[pdftex]{graphicx}
\else
  \usepackage[dvips]{graphicx}
\fi
\usepackage[cmex10]{amsmath}
\usepackage{algorithmic}
\usepackage[tight,footnotesize]{subfigure}
\usepackage{url}

\usepackage{psfrag}
\usepackage{amssymb}  
\usepackage[amsmath,thmmarks]{ntheorem}
\usepackage{multirow}
\usepackage{algorithm}
\usepackage{color,dsfont}

\theorembodyfont{\rmfamily} \theoremstyle{plain}
\newtheorem{theorem}{Theorem}

\newtheorem{corollary}{Corollary}
\newtheorem{definition}{Definition}
\newtheorem{proposition}{Proposition}
\newtheorem{assumption}{Assumption}
\newtheorem{remark}{Remark}
\newtheorem{example}{Example}
\newtheorem{problem}{Problem}

\newcommand{\subscr}[2]{#1_{\textup{#2}}}

\begin{document}
%
\title{\LARGE{Energy-Based Voronoi Partition\\in Constant Flow Environments}}
%
%
%

\author{Yu Ru and Sonia Martinez}
\maketitle

\begin{abstract}
As is well known, energy cost can greatly impact the deployment of battery-powered sensor networks in remote environments such as rivers or oceans. Motivated by this, we propose here an energy-based metric and associate energy-based Voronoi partitions with mobile vehicles in constant flows. The metric corresponds to the minimum energy that a vehicle requires to move from one point to another in the flow environment, and the resulting partition can be used by the vehicles in cooperative control tasks such as task assignment and coverage. Based on disk-based and asymptote-based approximations of the Voronoi regions, we determine a subset (or lower bound) and superset (or upper bound) of an agent's Voronoi neighbors. We then show that, via simulations, the upper bound is tight and its cardinality remains bounded as the number of generators increases. Finally, we propose efficient algorithms to compute the upper bound (especially when the generators dynamically change), which enables the fast calculation of Voronoi regions.
\end{abstract}


%
\IEEEpeerreviewmaketitle

\section{Introduction} \label{section1}
Due to the proliferation of low-cost sensing, communication, and computation devices, large groups of mobile vehicles equipped with sensors can be deployed into flow environments (e.g., rivers, lakes, oceans) to efficiently perform monitoring tasks. Depending on the nature of the application, multiple mobile vehicles can be coordinated based on different objectives. For search and rescue missions, a priority is to find/reach the target within the shortest time. However, for non-urgent tasks such as the monitoring of harmful algae blooms, maximizing the lifetime of the whole group of mobile vehicles can be more critical, as mobile vehicles are commonly powered by batteries with limited capacity. This motivates the study of minimum energy cooperative control algorithms for mobile vehicles in flow environments.

In this paper, we study a Voronoi partition associated with the minimum energy required for a vehicle to move from one point to another in a constant flow environment. We first derive an explicit expression for the energy-based metric, and study the Voronoi partition based on this metric using the vehicle locations as the set of generators. Similar to the time metric counterpart~\cite{yru_journal:Frazzoli_2004}, the Voronoi partition can then be used in the design of efficient target-assignment (or task allocation) algorithms (e.g., some other work~\cite{yru_journal:Sayyaadi_2011, yru_journal:Pavone_2011}). By assigning a vehicle to the targets that fall into its Voronoi region and guiding its motion appropriately, the vehicles can minimize the average energy spent by the group in servicing stochastic tasks that arrive according to a slow-rate Poisson distribution. However, contrary to the Euclidean case~\cite{yru_journal:Frazzoli_2004}, the Voronoi region defined by a general metric can be very involved. On the other hand, upper and lower approximations of the regions can be just enough to implement a coverage or target   assignment algorithm; see~\cite{yru_journal:Nowzari_2011}. Motivated by this, we propose methods to bound the set of Voronoi neighbors of a vehicle, which simplifies the calculation of Voronoi cells by vehicles. These are based on the following considerations: (i) the characterization of Voronoi region boundaries as hyperbolas, (ii) approximations of Voronoi regions by means of circles and polygons, and (iii) the derivation of a simple test (refer to Theorem~\ref{case_1}, Corollary~\ref{case_2}, and Proposition~\ref{prop:special_case}) that allows to discard vehicles that cannot be Voronoi neighbors. The test leads to an upper bound on the set of Voronoi neighbors of a vehicle. By generating vehicle locations (i.e., the set of generators for the Voronoi partition) independently according to a uniform distribution, we show that the average number of generators in the upper bound is bounded (via simulations) by $4.5$. Since the set of generators in the upper bound is sufficient for calculating Voronoi cells, the approach based on this upper bound (instead of using all generators) can save significant amount of time when calculating Voronoi cells, especially for applications with large amount of mobile vehicles. Therefore, we propose different algorithms to calculate the upper bound, especially when dealing with constantly changing generators.

Previously, Voronoi partitions in flow environments have been studied in connection to the shortest traveling time metric~\cite{yru_journal:Zermelo_1931, yru_journal:McGee_2006, yru_journal:Techy_2009, yru_journal:Kwok_2010, yru_journal:Bakolas_2010_Voronoi_j}. In contrast, there are relatively fewer works on the energy metric~\cite{yru_journal:Bongiorno_1967, yru_journal:Rowe_1990, yru_journal:Rowe_2000, yru_journal:Sun_2005, yru_journal:Ru_2011_b}. For example, in~\cite{yru_journal:Rowe_1990, yru_journal:Rowe_2000, yru_journal:Sun_2005}, the goal is to find a path with the minimum energy loss between given source and destination points in piecewise constant regions; in~\cite{yru_journal:Ru_2011_b}, the minimum energy metric (which is the same as this work) is used but the flow is modeled as a quadratic function (in this case, there is no explicit expression for the energy metric). In terms of approximating Voronoi cells, the work in~\cite{yru_journal:Evans_2008, yru_journal:Nowzari_2011} propose methods to deal with Voronoi partitions induced by the Euclidean distance metric (called standard Voronoi partitions). In~\cite{yru_journal:Cao_2003_rep, yru_journal:Bash_2007, yru_journal:Alsalih_2008}, distributed algorithms to calculate the standard Voronoi partitions are provided. For example, in~\cite{yru_journal:Cao_2003_rep}, explicit stopping criteria are proposed for a generator to calculate its own Voronoi cell without message broadcasting or routing. In contrast, the work in \cite{yru_journal:Bash_2007, yru_journal:Alsalih_2008} requires explicitly broadcasting generators or geographic routing. In our work, we assume that generator information is available (via either direct sensing or communication). As discussed in Section~\ref{section3_special}, the energy-based Voronoi partition can also be obtained via a slightly modified metric, which happens to be an additively weighted metric as in~\cite{yru_journal:Okabe_2000}. Methods such as~\cite{yru_journal:Fortune_1987, yru_journal:Karavelas_2002} can then be used to calculate the Voronoi partition based on the modified metric in a centralized fashion; however, to the best of our knowledge, there is no known distributed algorithm on calculating such partitions.

The contributions of this work are the following. i) An energy metric is proposed to study Voronoi partitions, which arises naturally in battery powered mobile vehicle applications in flow environments. To the best of our knowledge, this is the first work on Voronoi partition based on the minimum energy required for a vehicle to move from one point to another. In contrast, the traveling time based metric has been studied extensively for constant flows~\cite{yru_journal:Zermelo_1931, yru_journal:McGee_2006, yru_journal:Techy_2009}, piecewise constant flows~\cite{yru_journal:Kwok_2010}, and time varying flows~\cite{yru_journal:Bakolas_2010_Voronoi_j}. ii) In addition to deriving the lower bound on the set of Voronoi neighbors using a  disk-based lower approximation of Voronoi cells, an upper bound on the set of Voronoi neighbors is proposed utilizing asymptote-based lower and upper approximations of Voronoi cells. When deriving the upper bound, we introduce a dominance relation among Voronoi generators and provide a complete characterization for the dominance relation. iii) Since the upper bound is essential for a vehicle to compute its own Voronoi cell, we propose an efficient algorithm based on sorting generators (i.e., Algorithm~\ref{algorithm_efficient}) besides the method based on checking generators sequentially (i.e., Algorithm~\ref{algorithm_simple}). iv) To handle dynamically moving vehicles (namely, dynamically changing generators), we introduce a dominance graph for recomputing Voronoi cells only when absolutely necessary, which potentially avoids the recomputation due to any single change of the set of generators.

The paper is organized as follows. In Section~\ref{section2}, we define the minimum energy metric and formulate the Voronoi partition problem. Then we characterize the minimum energy metric and study the Voronoi partition for two generators in Section~\ref{section3}. In Section~\ref{section4}, we propose a disk based approximation for Voronoi cells and provide a lower bound on the set of Voronoi neighbors. To facilitate the calculation of Voronoi partitions in a distributed fashion, we study an asymptote based approximation of Voronoi cells, introduce the dominance relation and its characterization, and provide an upper bound on the set of Voronoi neighbors in Section~\ref{section5}. In Section~\ref{section6} we propose algorithms to calculate the upper bound, and show that the average number of generators in the upper bound is very small via simulations in Section~\ref{section7}. Finally, we summarize the work in Section~\ref{section8}.

\section{Problem Formulation} \label{section2}

In the Cartesian coordinate system, the studied flow environment is
described by $\mathds{R}^2$. The constant velocity field is a mapping $v: (x~y)^T \in \mathds{R}^2 \mapsto (B~0)^T$, where $B$ is a
positive constant. A vehicle runs at speed $U = (U_x~U_y)^T$ relative to the velocity field, and then the dynamic of the vehicle in the flow environment
can be described by
\begin{align}
\frac{dx}{dt} &= U_x + B~, \label{eq:dynamicsx}\\
\frac{dy}{dt} &= U_y~. \label{eq:dynamicsy}
\end{align}
We assume that vehicles can run against the flow.

To study Voronoi partitions, we introduce the following (pseudo)-metric.

\begin{definition}
Given two points $p^1$ and $p^2$ in the flow environment $\mathds{R}^2$, the energy metric $J(p^1, p^2)$ is defined as $J(p^1, p^2) = \min \int_{0}^{t_f} U^T U d t$, where $t_f$ is free, $U$ satisfies Eqs.~(\ref{eq:dynamicsx})
and (\ref{eq:dynamicsy}), and $x(0) = x_{p^1}$ (i.e., the $x$ coordinate of $p^1$), $y(0) = y_{p^1}$ (i.e., the $y$ coordinate of $p^1$), $x(t_f)
= x_{p^2}, ~y(t_f) = y_{p^2}$. \label{def:metric}
\end{definition}

The energy metric $J(p^1, p^2)$ is the minimum amount of energy required for the vehicle to move from its initial location $p^1$ to its final location $p^2$ among all possible controls. Note that there is no explicit constraint on $U$; however, as shown in Remark~\ref{remark:finite_energy}, the optimal control $U$ that achieves $J(p^1, p^2)$ is bounded by two times the flow velocity. The explicit expression for the energy metric is derived in Section~\ref{section3}. With this energy metric, now we can define the following Voronoi partition.

\begin{definition}
Let $P = \{p^1, p^2, ..., p^n\} \subset \mathds{R}^2$ be a set of distinct points, where $n \geq 2$. We call the region given by
$$V(p^i) = \{p \in \mathds{R}^2~|~ J(p^i, p) \leq J(p^j, p)~\mathrm{for}~j\neq i, j \in I_n \}$$
the energy-based Voronoi cell associated with $p^i$, where $I_n := \{1, 2, ..., n\}$, and the set given by $\mathcal{V} = \{V(p^1), V(p^2), ..., V(p^n)\}$ the energy-based Voronoi partition generated by $P$.
\label{def:Voronoi}
\end{definition}

Note that if $J(p^i, p)$ is replaced with $d_{p^i p} := \sqrt{(x_{p^i} - x_{p})^2 + (y_{p^i} - y_{p})^2}$, the partition is called a standard Voronoi partition. For simplicity, we use Voronoi partition to refer to the energy-based Voronoi partition in the rest of the paper.

Besides calculating the Voronoi partition $\mathcal{V}$ given the set of points $P$, we are especially interested in calculating each Voronoi cell $V(p^i)$ for $i = 1, ..., n$ given $\mathds{P} = P \setminus \{p^i\}$. For example, if $p^i$ is the location of a vehicle $V^i$ in the constant flow environment, $V(p^i)$ can be interpreted as the set of points that can be reached by $V^i$ with fewer energy consumption than by any other vehicle. If a task (e.g., taking measurements) has to be done at a point $p$ belonging to $V(p^i)$ and vehicle $V^i$ is assigned to the task, the energy consumption is minimized for this task. In this context, the challenge of computing the Voronoi cells lies in the fact that since the vehicles are moving, the generators $p^i$ for $i = 1, 2, ..., n$ are constantly changing. Without loss of generality, we formulate the following Voronoi cell calculation problem.

\begin{problem}
Given a fixed point $p^1$ and a set of points $\mathds{P} = \{p^2, p^3, ..., p^n\}$, calculate the Voronoi cell $V(p^1)$ as defined in Definition~\ref{def:Voronoi}. \label{problem:distributed_calculation}
\end{problem}

\section{Energy-Based Metric and Voronoi Partition with Two Generators} \label{section3}
In this section, we first study the minimum energy control problem and
provide an expression for the metric $J(p^1, p^2)$, and then derive the Voronoi boundary between two generators.

\subsection{Energy Metric: Expression}

The energy metric in Definition~\ref{def:metric} is given below.

\begin{proposition}
Given two points $p^1$ and $p^2$ in the flow environment $\mathds{R}^2$ with the velocity field $v$ satisfying $v_x(x, y) = B > 0$ and $v_y(x, y) = 0$, the
minimum energy $J(p^1, p^2) = \min \int_{0}^{t_f} U^T U d t$ is
\begin{equation}
J(p^1, p^2) = 2B(d_{p^1 p^2} + x_{p^1} - x_{p^2})~, \label{eq:min_energy}
\end{equation}
and the optimal control is $U(t) = -\frac{1}{2}
\begin{bmatrix} C_1\\C_2\end{bmatrix}$ for $t \in [0, t_f]$, where \begin{align}
C_1 = 2B(1 + \frac{x_{p^1} - x_{p^2}}{d_{p^1p^2}}),~~ C_2 = \frac{2B(y_{p^1} - y_{p^2})}{d_{p^1p^2}}, ~~t_f = \frac{d_{p^1 p^2}}{B}~.\label{eq:tf}
\end{align} \label{prop:metric}
\end{proposition}

\vspace{-30pt}\begin{proof}
Refer to the Appendix.
\end{proof}

\begin{remark}
Note that the quantity $J(p^1, p^2)$ is not a real metric because i) $J(p^1, p^2) = 0$ does not imply $p^1 = p^2$, and ii) $J(p^1, p^2)$ is not the same as $J(p^2, p^1)$ in general. More specifically, if $y_{p^1} = y_{p^2}$ and $x_{p^1} < x_{p^2}$, Eq.~\eqref{eq:min_energy} reduces to $2B(x_{p^2} - x_{p^1} + x_{p^1} - x_{p^2}) = 0$. This is consistent with the fact that no control is necessary if $p^2$ lies downstream of $p^1$. In addition, it can be verified that only when $x_{p^1} = x_{p^2}$, $J(p^1, p^2) = J(p^2, p^1)$. The magnitude of $U$ satisfies $\| U \| = \frac{\sqrt{C_1^2 + C_2^2}}{2} = B \sqrt{2 + \frac{2(x_{p^1} - x_{p^2})}{d_{p^1 p^2}}} \leq 2B$. Therefore, the optimal control is bounded. \hfill $\diamondsuit$ \label{remark:finite_energy}
\end{remark}



\subsection{Voronoi Partition: Two Generators} \label{section3_special}
Based on the metric $J(p^1, p^2)$ given in
Eq.~(\ref{eq:min_energy}), we now study the Voronoi partition given two generators $p^1, p^2$. Following Definition~\ref{def:Voronoi}, we have
$V(p^1) = \{p \in \mathds{R}^2~|~ J(p^1, p) \leq J(p^2, p)\}$. Using Eq.~\eqref{eq:min_energy}, $J(p^1, p) \leq J(p^2, p)$ can be rewritten as
\begin{align}
2B(d_{p^1 p} + x_{p^1} - x_{p}) &\leq 2B(d_{p^2 p} + x_{p^2} - x_{p})~, \nonumber\\
d_{p^1 p} + x_{p^1} &\leq d_{p^2 p} + x_{p^2}~, \label{eq:1}\\
d_{p^1 p} - d_{p^2 p} &\leq x_{p^2} - x_{p^1}~, \label{eq:3}
\end{align}
where we obtain Eq.~(\ref{eq:1}) because $B > 0$.

Based on Eq.~(\ref{eq:1}), we can also use the metric $d_{p^1 p} +
x_{p^1}$ to obtain the same Voronoi partition.
In~\cite{yru_journal:Okabe_2000}, this metric falls into the
category of additively weighted distances. Therefore, given the set of points $P$, the Voronoi partition in Definition~\ref{def:Voronoi} can be calculated using existing methods such as~\cite{yru_journal:Fortune_1987, yru_journal:Karavelas_2002}. Straightforward application of such methods to Problem~\ref{problem:distributed_calculation} can be very inefficient because all Voronoi cells have to be computed in order to just obtain $V(p^1)$. Another simple idea to solve Problem~\ref{problem:distributed_calculation} is that we consider one point at a time and keep refining $V(p^1)$ until all points have been taken into account. However, as we show in Section~\ref{section5}, only a subset of points are necessary for calculating $V(p^1)$. More details on comparing different methods to solve Problem~\ref{problem:distributed_calculation} are provided in Section~\ref{section7}.

Now we can rewrite $V(p^1)$ as $V(p^1) = \{p \in \mathds{R}^2~|~ d_{p^1 p} - d_{p^2 p} \leq x_{p^2} - x_{p^1}\}$, and the boundary between $p^1$ and $p^2$ is
$B(p^1, p^2) := \{p \in \mathds{R}^2~|~ d_{p^1 p} - d_{p^2 p} = x_{p^2} - x_{p^1}\}$. Without loss of generality, we assume that $x_{p^1} \leq x_{p^2}$ and $y_{p^1} \leq y_{p^2}$. Depending on the relative position between $p^1$ and $p^2$, there are three different cases for the Voronoi cells and the
boundary between these cells.

\begin{figure}[tb] \centering

\subfigure[Original coordinates.]{\psfrag{0}{$0$} \psfrag{3}{$p^1$} \psfrag{4}{$p^2$}
\psfrag{5}{$\alpha$} \psfrag{x}{$x$} \psfrag{y}{$y$} \psfrag{6}{$p^3$}
\includegraphics[scale=0.25]{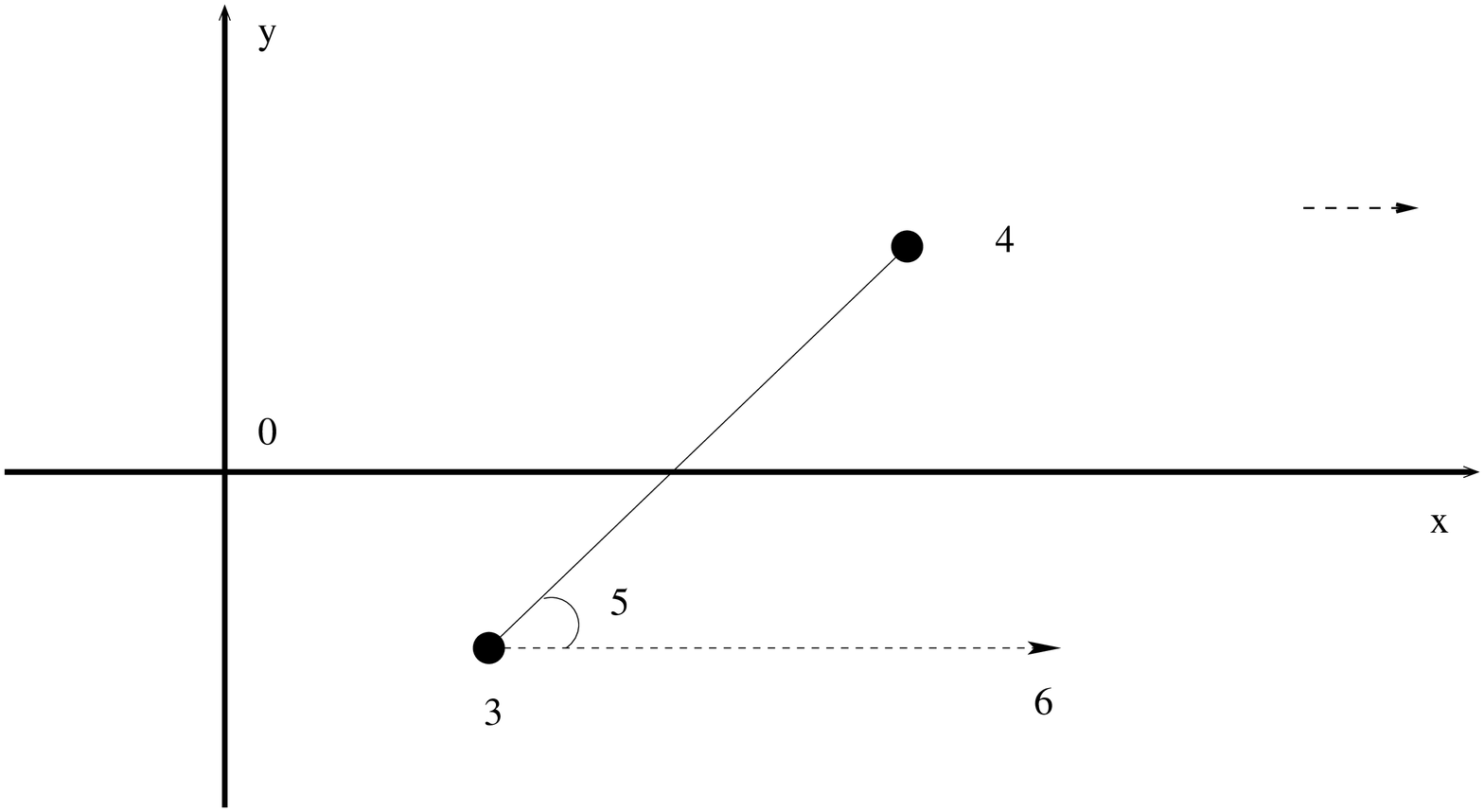}} \hfil \subfigure[Transformed coordinates.]{\psfrag{0}{$0$} \psfrag{1}{$y'$}
\psfrag{2}{$x'$} \psfrag{3}{$p^1$} \psfrag{4}{$p^2$}
\psfrag{5}{$p^*$} \psfrag{6}{$r_{p^1}$}
\psfrag{7}{$r_{p^2}$} \psfrag{8}{$l_1$} \psfrag{9}{$l_2$}
\psfrag{10}{$l_3$} \psfrag{11}{$l_4$}
\includegraphics[scale=0.25]{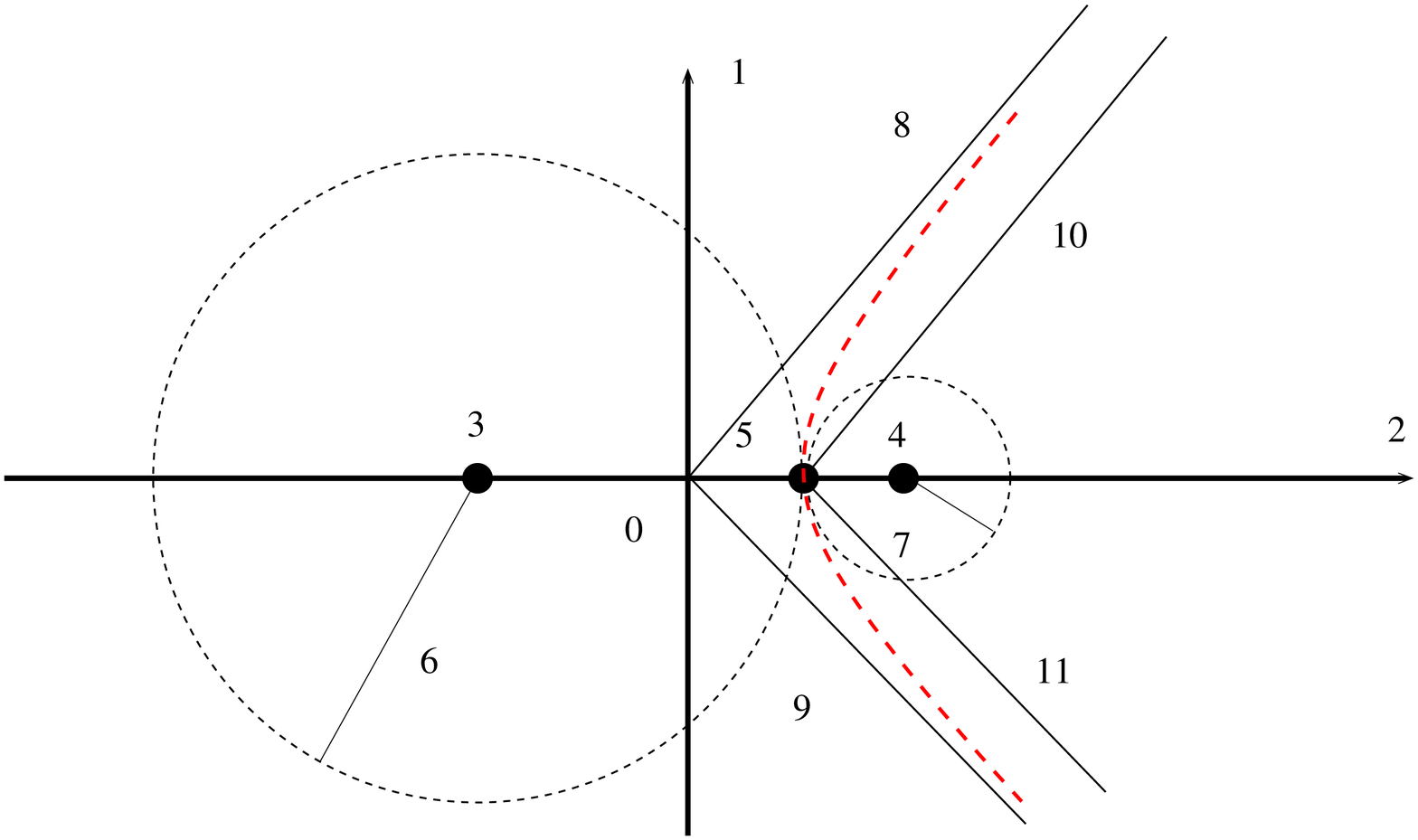}}

\caption{Voronoi partition generated by $\{p^1, p^2\}$ satisfying
$x_{p^1} < x_{p^2}$ and $y_{p^1} < y_{p^2}$. The red dashed line is
the boundary.} \label{fig0} \vspace{-15pt}
\end{figure}

\textbf{Case I}: $x_{p^1} < x_{p^2}$ \textbf{and} $y_{p^1} < y_{p^2}$. For any point $p$ on the boundary we have $d_{p^1 p} - d_{p^2 p} =
x_{p^2} - x_{p^1} > 0$. It can be verified that the boundary is a
hyperbolic curve. To derive an equation, we first transform the
coordinate from $(x, y)$ to $(x', y')$ such that the origin is at
$(\frac{x_{p^1} + x_{p^2}}{2}, \frac{y_{p^1} + y_{p^2}}{2})$ and the
positive $x'$ direction is from $p^1$ to $p^2$. Essentially the
transformation involves shifting the origin and rotating the $x, y$ axes. As shown in Fig.~\ref{fig0}(a), we use $\alpha$ to denote the
rotation angle $\angle p^3p^1p^2$, where $p^1p^3$ is parallel to the $x$ axis, and we have
$\tan \alpha = \frac{y_{p^2} - y_{p^1}}{x_{p^2} - x_{p^1}}$. Then
for any point $p$ with coordinates $(x_p, y_p)$, its coordinate in
the $(x', y')$ plane is given as
\begin{align}
x_p' &= (x_p - \frac{x_{p^1} + x_{p^2}}{2}) \cos \alpha + (y_p - \frac{y_{p^1} + y_{p^2}}{2}) \sin \alpha~, \label{eq:x_new}\\
y_p' &= -(x_p - \frac{x_{p^1} + x_{p^2}}{2}) \sin \alpha + (y_p -
\frac{y_{p^1} + y_{p^2}}{2}) \cos \alpha~. \label{eq:y_new}
\end{align}
In the transformed coordinate, the boundary is shown as the red dotted hyperbola in Fig.~\ref{fig0}(b), and can be described by the equation
$\frac{(x_p')^2}{a^2} - \frac{(y_p')^2}{b^2} = 1$, where $a =
\frac{x_{p^2} - x_{p^1}}{2}$, $c = \frac{d_{p^1 p^2}}{2}$, $b = \sqrt{c^2
- a^2} = \frac{y_{p^2} - y_{p^1}}{2}$, and $x_p' \geq a$. Note that
in Fig.~\ref{fig0}(b), the coordinate for $p^*$ (namely, the
intersection point between the boundary and the $x'$ axis) is $(a,
0)$. Therefore, in the original coordinate, the boundary can be
described as
\begin{equation} \resizebox{.9\hsize}{!}{$\frac{((x_p - \frac{x_{p^1} + x_{p^2}}{2}) \cos \alpha + (y_p -
\frac{y_{p^1} + y_{p^2}}{2}) \sin \alpha)^2}{(\frac{x_{p^2} -
x_{p^1}}{2})^2} - \frac{(-(x_p - \frac{x_{p^1} + x_{p^2}}{2}) \sin
\alpha + (y_p - \frac{y_{p^1} + y_{p^2}}{2}) \cos
\alpha)^2}{(\frac{y_{p^2} - y_{p^1}}{2})^2} = 1$}\label{eq:boundary}
\end{equation} with the constraint that
\begin{equation} (x_p - \frac{x_{p^1} + x_{p^2}}{2}) \cos
\alpha + (y_p - \frac{y_{p^1} + y_{p^2}}{2}) \sin \alpha \geq \frac{x_{p^2} - x_{p^1}}{2}~.
\label{eq:boundary_constraint}
\end{equation}

Now we apply the above equations of the boundary to specific scenarios.

\begin{figure}[tb] \centering

\subfigure[$x_{p^1} = x_{p^2}$ and $y_{p^1} < y_{p^2}$]{\psfrag{0}{$0$}
\psfrag{3}{$p^2$} \psfrag{4}{$p^1$} \psfrag{x}{$x$} \psfrag{y}{$y$}
\includegraphics[scale=0.25]{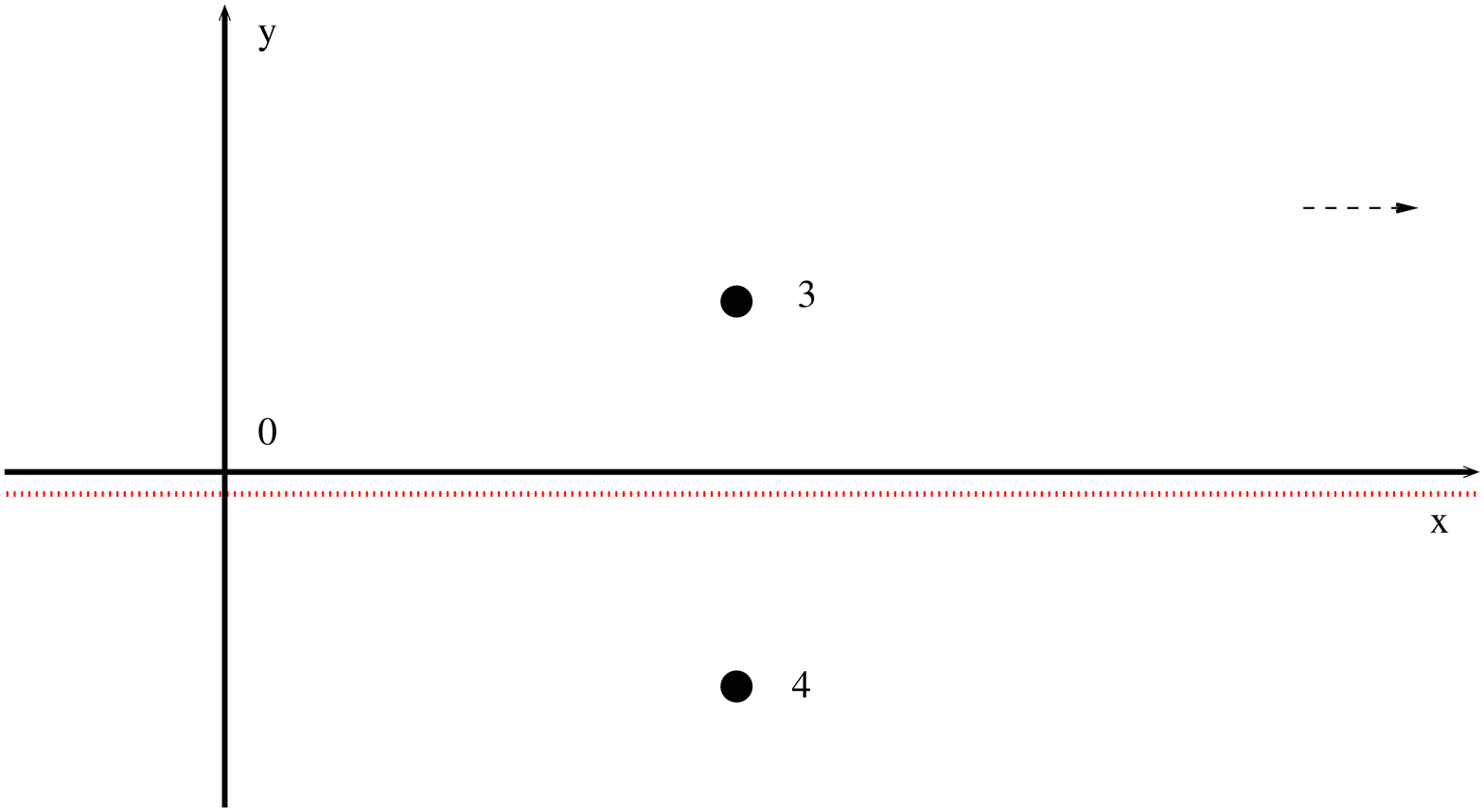}} \hfil \subfigure[$x_{p^1} < x_{p^2}$ and $y_{p^1} =
y_{p^2}$.]{\psfrag{0}{$0$}
\psfrag{3}{$p^1$} \psfrag{4}{$p^2$} \psfrag{x}{$x$} \psfrag{y}{$y$}
\includegraphics[scale=0.25]{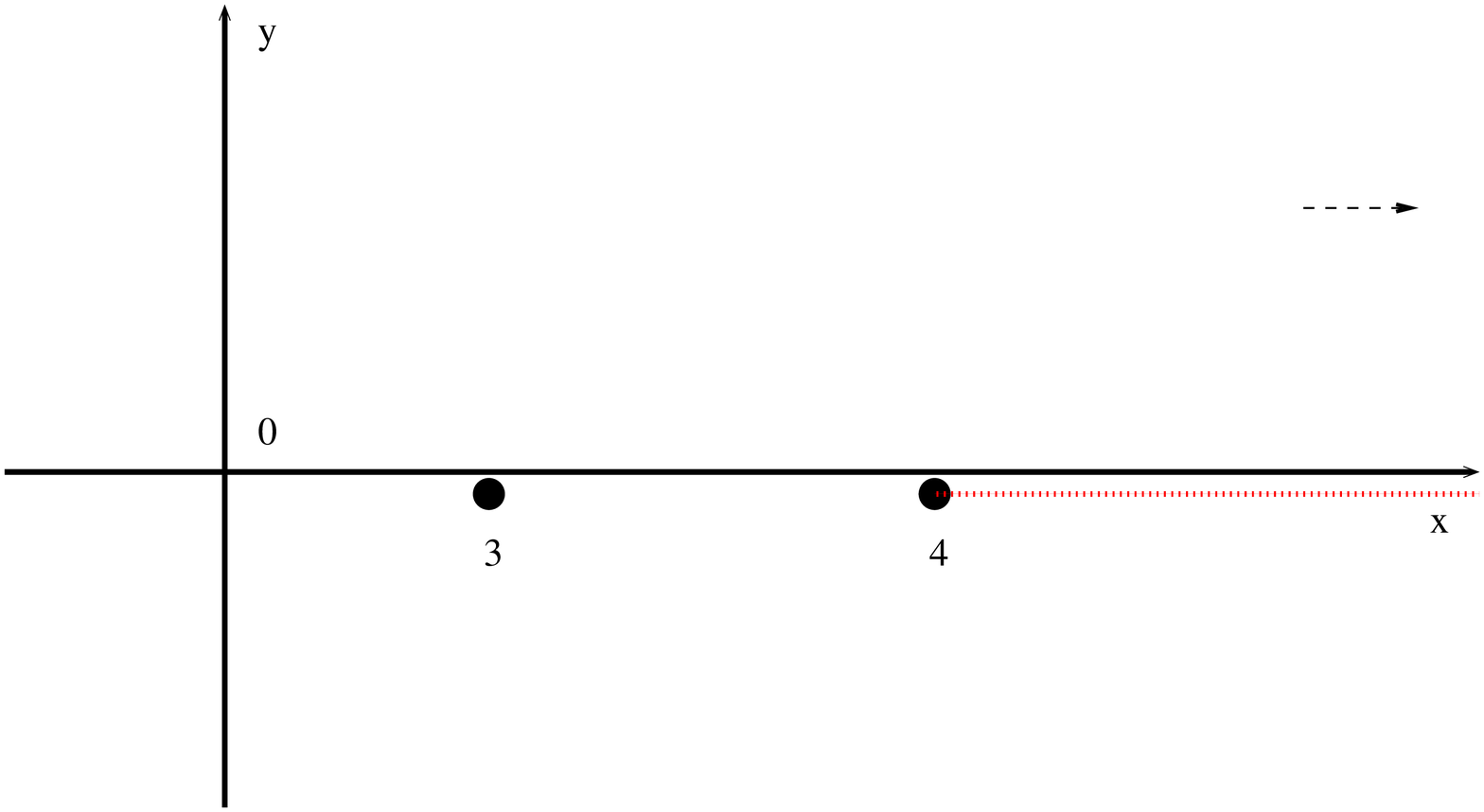}}
\caption{Voronoi partition generated by $\{p^1, p^2\}$. The red
dashed line is the boundary.} \label{fig1} \vspace{-15pt}
\end{figure}

\textbf{Case II}: $x_{p^1} = x_{p^2}$ \textbf{and} $y_{p^1} < y_{p^2}$. In this case,  $\alpha = \frac{\pi}{2}$, Eq.~(\ref{eq:boundary})
reduces to $y_p = \frac{y_{p^1} + y_{p^2}}{2}$, and
Eq.~(\ref{eq:boundary_constraint}) holds trivially. In other words,
the boundary is a perpendicular bisector of the line segment
$p^1p^2$, as shown in Fig.~\ref{fig1}(a).

\textbf{Case III}: $x_{p^1} < x_{p^2}$ \textbf{and} $y_{p^1} = y_{p^2}$. In this case,  $\alpha = 0$, Eq.~(\ref{eq:boundary}) reduces to
$y_p = \frac{y_{p^1} + y_{p^2}}{2} = y_{p^2}$, and
Eq.~(\ref{eq:boundary_constraint}) reduces to $x_{p} \geq x_{p^2}$.
The boundary is a half line given by $\{p \in \mathds{R}^2~|~x_{p} \geq
x_{p^2}, y_{p} = y_{p^2}\}$, as shown in Fig.~\ref{fig1}(b).

It is straightforward to obtain similar equations for the cases with $x_{p^1} > x_{p^2}$ and/or $y_{p^1} > y_{p^2}$.

\section{Disk-Based Lower Approximation of Voronoi Cells and Lower Bound on Voronoi Neighbors} \label{section4}
In this section, we first study disk-based lower approximation of Voronoi cells, and then derive a lower bound on the set of Voronoi neighbors.

\subsection{Disk-Based Lower Approximation}
The disk-based lower approximation of $V(p^i)$ is given as $D(p^i, r_{p^i}) = \{p \in \mathds{R}^2~|~d_{p^i p} \leq r_{p^i}\}$, i.e., a disk centered at $p^i$ with radius $r_{p^i}$, such that $D(p^i, r_{p^i}) \subseteq V(p^i)$.

We first study the case with two points $p^1$ and $p^2$ satisfying $x_{p^1} \leq x_{p^2}$. Since in general the boundary between the two Voronoi cells is a hyperbola as shown in Fig.~\ref{fig0}(b), the radius for $p^1$ can be chosen to be $r_{p^1} = c + a = \frac{d_{p^1 p^2}}{2} + \frac{x_{p^2} - x_{p^1}}{2}$, and the radius for $p^2$ can be chosen to be $r_{p^2} = c - a = \frac{d_{p^1 p^2}}{2} + \frac{x_{p^1} - x_{p^2}}{2}$. Note that the hyperbola boundary intersects with the disk $D(p^2, r_{p^2})$ only at one point (namely, the point $p^*$ in Fig.~\ref{fig0}(b); we denote the point as $p^*(p^1, p^2)$) because the focus of the hyperbola is the same as the center of the disk and the eccentricity of a hyperbola is larger than $1$. We have $r_{p^2} = d_{p^2 p^*(p^1, p^2)}$. It can be verified that the hyperbola boundary intersects with the disk $D(p^1, r_{p^1})$ also only at the point $p^*(p^1, p^2)$ and $r_{p^1} = d_{p^1 p^*(p^1, p^2)}$. In the original coordinate, we have $p^*(p^1, p^2) = p^1 + \frac{c+a}{2c} (p^2 - p^1)$, where $a = \frac{x_{p^2} - x_{p^1}}{2}$ and $c = \frac{d_{p^1 p^2}}{2}$.

If $x_{p^1} = x_{p^2}$ and $y_{p^1} < y_{p^2}$ (i.e., \textbf{Case II} as discussed in Section~\ref{section3_special}), then we have $r_{p^1} = r_{p^2} = \frac{d_{p^1 p^2}}{2}$. If $x_{p^1} < x_{p^2}$ and $y_{p^1} = y_{p^2}$ (i.e., \textbf{Case III} as discussed in Section~\ref{section3_special}), then we have $r_{p^1} = d_{p^1 p^2}$ and $r_{p^2} = 0$. It can be verified that the above results also hold if $x_{p^1} > x_{p^2}$.

If there are $n$ points $p^1, p^2, ..., p^n$, we can choose the radius $r_{p^i} = \min_{j \in \{1, 2, ..., n\} \setminus \{i\}} (\frac{d_{p^i p^j}}{2} + \frac{x_{p^j} - x_{p^i}}{2})$. Therefore, $D(p^i, r_{p^i}) \subseteq V(p^i)$.

\subsection{Lower Bound on Voronoi Neighbors}
Now we focus on $V(p^1)$ and let
\begin{equation}
\mathcal{N}_D(p^1) = \underset{p^j \in \{p^2, ..., p^n\}}{\operatorname{argmin}} (\frac{d_{p^1 p^j}}{2} + \frac{x_{p^j} - x_{p^1}}{2})~. \label{eq:NDp1}
\end{equation} In general, $\mathcal{N}_D(p^1)$ could be a set with multiple elements. Before analyzing $\mathcal{N}_D(p^1)$, we first define the set of Voronoi neighbors of $p^1$.

\begin{definition}
Given a fixed point $p^1$ and a set of points $\mathds{P} = \{p^2, p^3, ..., p^n\}$, a point $p \in \mathds{P}$ is a Voronoi neighbor of $p^1$ if $V(p^1) \cap V(p)$ is non-empty and non-trivial (i.e., not a single point). We use $\mathcal{N}_{V}(p^1)$ to denote the set of Voronoi neighbors of $p^1$. \label{def:Voronoi_neighbor}
\end{definition}

\begin{remark}
Note that in our setting, $V(p^1) \cap V(p)$ could be empty, a single point, a line segment, or part of a hyperbola. In our definition of Voronoi neighbors, we treat two points as neighbors only when the intersection of their Voronoi cells is nonempty and nontrivial (i.e., not a single point). This is because we are primarily interested in calculating Voronoi cells and ruling out this trivial case does not affect the calculation. A definition of Voronoi neighbors in the same spirit is used in~\cite{yru_journal:Bakolas_2010_Voronoi_j}. \hfill $\diamondsuit$
\end{remark}

Now we are ready to state the relationship between $\mathcal{N}_D(p^1)$ and $\mathcal{N}_V(p^1)$.

\begin{theorem}
Given a fixed point $p^1$ and a set of points $\mathds{P} = \{p^2, p^3, ..., p^n\}$, $\mathcal{N}_D(p^1) \subseteq \mathcal{N}_V(p^1)$, where $\mathcal{N}_D(p^1)$ is defined in Eq.~\eqref{eq:NDp1}. \label{theorem:Voronoi_neighbor_lower}
\end{theorem}

\vspace{-10pt}\begin{proof}
Refer to the Appendix.
\end{proof}

Now it is clear that $\mathcal{N}_D(p^1)$ is a lower bound on the set of Voronoi neighbors of $p^1$. Even if a point $p^k$ may not minimize the radius of the disk that lower approximates $V(p^1)$, $p^k$ can still be a Voronoi neighbor of $p^1$. Thus $\mathcal{N}_D(p^1)$ could potentially be augmented to a larger set $\overline{\mathcal{N}_D(p^1)}$ while still satisfying $\overline{\mathcal{N}_D(p^1)} \subseteq \mathcal{N}_V(p^1)$. It can be verified that the following sufficient condition for testing if $p^k$ is a Voronoi neighbor of $p^1$ holds.

\begin{proposition}
Given a fixed point $p^1$ and a set of points $\mathds{P} = \{p^2, p^3, ..., p^n\}$, a point $p^k$ for $k \in \{2, 3, ..., n\}$ is a Voronoi neighbor of $p^1$ if for any $l \in \{2, 3, ..., n\} \setminus \{k\}$ we have $J(p^k, p^*(p^1, p^k)) < J(p^l, p^*(p^1, p^k))$. \label{prop:neighbor_lower}
\end{proposition}

\begin{proof}
Refer to the Appendix.
\end{proof}

Essentially, the condition verifies that for a point $p^k$, the energy required to reach the specific point $p^*(p^1, p^k)$ (which belongs to the boundary between $p^1$ and $p^k$) from any point other than $p^k$ is strictly larger than the energy from $p^k$. Based on this condition, we can construct $\overline{\mathcal{N}_D(p^1)}$ by starting with $\mathcal{N}_D(p^1)$, and adding a point $p^k$ to $\mathcal{N}_D(p^1)$ if the condition in Proposition~\ref{prop:neighbor_lower} holds.



\section{Asymptote-Based Approximations of Voronoi Cells and Upper Bound on Voronoi Neighbors} \label{section5}

In this section, we first propose asymptote-based lower and upper approximations of Voronoi cells, and then introduce a dominance relation to upper bound the set of Voronoi neighbors.

\subsection{Asymptote-Based Approximations of Voronoi Cells}
Since in general the boundary is a hyperbola, another way to approximate the Voronoi cells is to use the asymptotes of the hyperbola. Here, we focus on two points $p^1$ and $p^2$ satisfying $x_{p^1} < x_{p^2}$ and $y_{p^1} < y_{p^2}$. In the $x'-y'$ plane of Fig.~\ref{fig0}(b), the equation for the asymptote $l_1$ (or $l_2$) is $y' = \frac{b}{a}x'$ (or $y' = -\frac{b}{a}x'$), where $a = \frac{x_{p^2} - x_{p^1}}{2}$ and $b = \frac{y_{p^2} - y_{p^1}}{2}$. It can be shown that the region described by $\subscr{D}{lower}'(p^1 | p^2) = \{(x'~y')^T \in \mathds{R}^2~|~x' \leq \frac{a}{b} y'~\mathrm{if}~y' \geq 0, x' \leq -\frac{a}{b} y'~\mathrm{if}~y' < 0\}$ satisfies $\subscr{D}{lower}'(p^1 | p^2) \subseteq V'(p^1)$, where $V'(p^1)$ is the Voronoi cell $V(p^1)$ in the transformed $x'-y'$ plane. Going back to the original $x-y$ plane, we have $\subscr{D}{lower}(p^1 | p^2) \subseteq V(p^1)$. At the same time, we get an upper approximation for $V(p^2)$ as $\subscr{D}{upper}(p^2 | p^1) = \mathds{R}^2 \setminus \subscr{D}{lower}(p^1 | p^2)$.

To obtain an upper approximation for $V(p^1)$, we use $l_3$ (which is parallel to $l_1$) and $l_4$ (which is parallel to $l_2$) that pass through the point $p^*$ in Fig.~\ref{fig0}(b). The equation for $l_3$ (or $l_4$) is $y' = \frac{b}{a} (x' - a)$ (or $y' = -\frac{b}{a} (x' - a)$). It can be shown that the region described by $\subscr{D}{upper}'(p^1|p^2) = \{(x'~y')^T \in \mathds{R}^2~|~x' \leq \frac{a}{b} y' + a~\mathrm{if}~y' \geq 0, x' \leq -\frac{a}{b} y' + a~\mathrm{if}~y' < 0\}$ satisfies $V'(p^1) \subseteq \subscr{D}{upper}'(p^1|p^2)$. Going back to the original $x-y$ plane, we have $V(p^1) \subseteq \subscr{D}{upper}(p^1|p^2)$. At the same time, we get a lower approximation for $V(p^2)$ as $\subscr{D}{lower}(p^2) = \mathds{R}^2 \setminus \subscr{D}{upper}(p^1)$.

If there are $n$ points, we can lower approximate $V(p^i)$ using $\subscr{V}{lower}(p^i) = \cap_{j \in I_n \setminus \{i\}} \subscr{D}{lower}(p^i | p^j)$, and upper approximate $V(p^i)$ using $\subscr{V}{upper}(p^i) = \cup_{j \in I_n \setminus \{i\}} \subscr{D}{upper}(p^i | p^j)$, where $D(p^i | p^j)$ is the approximation of $V(p^i)$ given the point $p^j$.

Since $l_1, \dots, l_4$ play a very important role in the approximations, we study their equations in the $x-y$ plane. Here we focus on $l_1$ and $l_2$ since $l_3$ (or $l_4$) is parallel to $l_1$ (or $l_2$). For $l_1$, we are interested in $y' = \frac{b}{a}x'$ for $y' \geq 0$, while for $l_2$, we are interested in $y' = -\frac{b}{a}x'$ for $y' \leq 0$.

\begin{proposition}
Given two points $p^1$ and $p^2$ in the flow environment $\mathds{R}^2$ satisfying $x_{p^1} < x_{p^2}$ and $y_{p^1} < y_{p^2}$, two asymptotes of the hyperbolic boundary between the Voronoi cells of $p^1$ and $p^2$ are
\begin{equation}
y - \frac{y_{p^1} + y_{p^2}}{2} = (x - \frac{x_{p^1} + x_{p^2}}{2}) \tan2 \alpha~\mathrm{with}~y \geq \frac{y_{p^1} + y_{p^2}}{2}~, \label{eq:l1}
\end{equation}
where $\alpha = \arctan  \frac{y_{p^2} - y_{p^1}}{x_{p^2} - x_{p^1}}$,
and
\begin{equation}
y = \frac{y_{p^1} + y_{p^2}}{2}~\mathrm{and}~x \geq \frac{x_{p^1} + x_{p^2}}{2}~. \label{eq:l2}
\end{equation} \label{prop:asymptotes}
\end{proposition}

\vspace{-30pt}\begin{proof}
Refer to the Appendix.
\end{proof}

\begin{remark}
If $x_{p^1} < x_{p^2}$ and $y_{p^1} = y_{p^2}$, which implies that $\alpha = 0$, then Eq.~\eqref{eq:l1} becomes $y = \frac{y_{p^1} + y_{p^2}}{2}$, and $x \geq \frac{x_{p^1} + x_{p^2}}{2}$ because $x' \geq 0$, while Eq.~\eqref{eq:l2} remains the same. In this case, the two asymptotes coincide and form the exact boundary as shown in Fig.~\ref{fig1}(b). If $x_{p^1} = x_{p^2}$ and $y_{p^1} < y_{p^2}$, which implies that $\alpha = \frac{\pi}{2}$, then Eq.~\eqref{eq:l1} becomes $y = \frac{y_{p^1} + y_{p^2}}{2}$, and $x \leq \frac{x_{p^1} + x_{p^2}}{2}$ because $y' \geq 0$, while Eq.~\eqref{eq:l2} remains the same. In this case, the two asymptotes form a straight line and are also the exact boundary as shown in Fig.~\ref{fig1}(a). \hfill $\diamondsuit$
\end{remark}

\begin{remark}
Note that the two asymptotes pass through the middle point of $p^1p^2$. One is always parallel to the $x$ axis, while the other has the slope $\tan 2 \alpha$. \hfill $\diamondsuit$
\end{remark}


Similarly, we can obtain asymptote equations if $x_{p^1} \geq x_{p^2}$ and/or $y_{p^1} \geq y_{p^2}$. In the next subsection, we introduce a dominance relation and provide conditions to check the dominance relation, in which the asymptote equations prove to be useful.

\subsection{Dominance Relation and its Characterization} \label{section5_scenario}

When we calculate the Voronoi cell $V(p^1)$ in Problem~\ref{problem:distributed_calculation} by considering point $p^2$ first and then $p^3$, it is possible that the Voronoi cell of $p^1$ is not strictly refined when considering $p^3$ given $p^2$; in this scenario, $p^3$ is not necessary to compute $V(p^1)$, and potentially the calculation of $V(p^1)$ can be done more efficiently. Now we introduce a dominance relation which captures this scenario.

\begin{definition}
Given a fixed point $p^1$, and two points $p^2, p^3$ (that are different from $p^1$), we say $p^2$ dominates $p^3$ (denoted as $p^2 \succ p^3$) if $V(p^1~|~p^2) = V(p^1~|~p^2, p^3)$, where $V(p^1~|~p^2) = \{p \in \mathds{R}^2~|~J(p^1, p) \leq J(p^2, p)\}$ and $V(p^1~|~p^2, p^3) = \{ p \in \mathds{R}^2~|~J(p^1, p) \leq J(p^2, p), J(p^1, p) \leq J(p^3, p)\}$. \label{def:dominance}
\end{definition}

By definition, $p^2$ dominates itself; if $p^2$ dominates $p^3$, then only $p^2$ matters when $p^1$ calculates its Voronoi cell (this is proved more generally in Theorem~\ref{prop:Voronoi_neighbor}). Given a fixed point $p^1$, to check if $p^2$ dominates $p^3$, there are four scenarios:
\begin{itemize}
\item \textbf{Scenario A}: $y_{p^1} < y_{p^2}$. The iff condition is given in Theorem~\ref{case_1};
\item \textbf{Scenario B}: $y_{p^1} > y_{p^2}$. The iff condition is given in Corollary~\ref{case_2};
\item \textbf{Scenario C}: $y_{p^1} = y_{p^2}$ and $x_{p^1} > x_{p^2}$. The iff condition is given in part (a) of Proposition~\ref{prop:special_case};
\item \textbf{Scenario D}:  $y_{p^1} = y_{p^2}$ and $x_{p^1} < x_{p^2}$. The iff condition is given in part (b) of Proposition~\ref{prop:special_case}.
\end{itemize}

\begin{theorem}
Given a fixed point $p^1$, and two points $p^2, p^3$ that are different from $p^1$ and satisfy $y_{p^1} < y_{p^2}$, $p^2 \succ p^3$ iff
$y_{p^3} \geq y_{p^2}$, and
\begin{equation}
(y_{p^3} - y_{p^1}) \times (x_{p^2} - x_{p^1}) \leq (x_{p^3} - x_{p^1}) \times (y_{p^2} - y_{p^1})~.\label{eq:inequality}
\end{equation} \label{case_1}
\end{theorem}

\vspace{-30pt}\begin{proof}
Refer to the Appendix.
\end{proof}

Similarly, we can prove the following result for the case $y_{p^1} > y_{p^2}$.

\begin{corollary}
Given a fixed point $p^1$, and two points $p^2, p^3$ that are different from $p^1$ and satisfy $y_{p^1} > y_{p^2}$, $p^2 \succ p^3$ iff
$y_{p^3} \leq y_{p^2}$, and $(y_{p^3} - y_{p^1}) \times (x_{p^2} - x_{p^1}) \geq (x_{p^3} - x_{p^1}) \times (y_{p^2} - y_{p^1})$. \label{case_2}
\end{corollary}

If $y_{p^1} = y_{p^2}$ for points $p^1$ and $p^2$, the following result can be verified.
\begin{proposition}
Given a fixed point $p^1$, and two points $p^2, p^3$ that are different from $p^1$ and satisfy $y_{p^1} = y_{p^2}$ and
\begin{itemize}
\item[(a)] $x_{p^1} > x_{p^2}$, $p^2 \succ p^3$ iff $y_{p^3} \neq y_{p^2}$, or $y_{p^3} = y_{p^2}$ and $x_{p^3} < x_{p^1}$.
\item[(b)] $x_{p^1} < x_{p^2}$, $p^2 \succ p^3$  iff
$y_{p^3} = y_{p^2}$ and $x_{p^3} \geq x_{p^2}$.
\end{itemize} \label{prop:special_case}
\end{proposition}

\begin{figure}[tb] \centering
\psfrag{3}{$p^1$} \psfrag{4}{$p_a^2$}
\psfrag{5}{$p_b^2$} \psfrag{6}{$p_c^2$}
\includegraphics[scale=0.35]{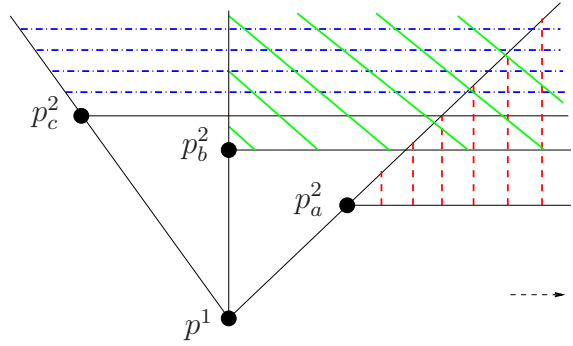}
\caption{Three different cases for Scenario A ($y_{p^1} < y_{p^2}$): $x_{p_a^2} > x_{p^1}$, $x_{p_b^2} = x_{p^1}$, and $x_{p_c^2} < x_{p^1}$.} \label{fig6} \vspace{-15pt}
\end{figure}

\begin{remark}
Fig.~\ref{fig6} illustrates the region of point $p^3$ that is dominated by $p^2$ in Scenario A. There are three different cases depending on the $x$ coordinates of $p^1$ and $p^2$. If $x_{p^2} > x_{p^1}$ (or $x_{p^2} = x_{p^1}$, $x_{p^2} < x_{p^1}$, respectively), any point in the red dotted (or green solid, blue dashed, respectively) region in Fig.~\ref{fig6} is dominated by $p^2$. \hfill $\diamondsuit$
\end{remark}


%

\subsection{Upper Bound on Voronoi Neighbors}

In this subsection, we propose an upper bound on the set of Voronoi neighbors based on the dominance relation introduced earlier. We first show that the dominance relation is antisymmetric under certain conditions.

\begin{proposition} (Antisymmetry of Dominance)
Given a fixed point $p^1$, and two points $p^2, p^3$ that are different from $p^1$, if $y_{p^1} \neq y_{p^2}$, or $y_{p^1} = y_{p^2}$ and $x_{p^1} < x_{p^2}$, then $p^2 \succ p^3$ and $p^3 \succ p^2$ implies that $p^2 = p^3$. \label{prop:antisymmetry}
\end{proposition}

\begin{proof}
Refer to the Appendix.
\end{proof}

Given a fixed point $p^1$, and two points $p^2, p^3$ that are different from $p^1$, if $y_{p^1} = y_{p^2}$ and $x_{p^1} > x_{p^2}$ (namely, Scenario C in Section~\ref{section5_scenario}), $p^2 \succ p^3$ and $p^3 \succ p^2$ may not imply $p^2 = p^3$. In fact, as long as $y_{p^1} = y_{p^2} = y_{p^3}$, $x_{p^1} > x_{p^2}$ and $x_{p^1} > x_{p^3}$, we have $p^2 \succ p^3$ and $p^3 \succ p^2$. The reason is that for any $p$ satisfying $y_{p} = y_{p^1}$ and $x_{p} < x_{p^1}$, $V(p^1|p) = \{(x~y)^T \in \mathds{R}^2~|~y = y_{p^1}, x\geq x_{p^1}\}$ which does not rely on the exact location of $p$. In this case, the Voronoi cell of $p^1$ is degenerated. To simplify the discussion, we make the following assumption, which guarantees that the Voronoi cell of $p^1$ is nonempty.

\begin{assumption}
Given a fixed point $p^1$, and a set of points $\mathds{P} = \{p^2, p^3, ..., p^{n}\}$ satisfying points in $\mathds{P} \cup \{p^1\}$ are distinct, $\forall p \in \mathds{P}$, it holds that $y_{p} \neq y_{p^1}$, or $y_{p} = y_{p^1}$ and $x_{p} > x_{p^1}$. \label{assumption}
\end{assumption}

Given a fixed point $p^1$ and two points $p^2, p^3$, if $p^1, p^2, p^3$ satisfy Assumption~\ref{assumption}, Proposition~\ref{prop:antisymmetry} guarantees that there are only three cases in terms of the dominance relation between $p^2$ and $p^3$, i.e., $p^2 \succ p^3$, $p^3 \succ p^2$, or $p^2$ and $p^3$ do not dominate each other.

\begin{definition}
Given a fixed point $p^1$ and a set of points $\mathds{P} = \{p^2, p^3, ..., p^n\}$ satisfying Assumption~\ref{assumption}, we define $\mathcal{N}_{G}(p^1)$ as the set of points satisfying that for any $p \in \mathcal{N}_{G}(p^1)$  there does not exist another point $p' \in \mathds{P}$ that is different from $p$ and dominates $p$. \label{def:NGp1}
\end{definition}

Now we are ready to state the relationship between $\mathcal{N}_G(p^1)$ and $\mathcal{N}_V(p^1)$.

\begin{theorem}
Given a fixed point $p^1$ and a set of points $\mathds{P} = \{p^2, p^3, ..., p^n\}$ satisfying Assumption~\ref{assumption}, $\mathcal{N}_V(p^1) \subseteq \mathcal{N}_G(p^1)$. \label{prop:Voronoi_neighbor}
\end{theorem}

\begin{proof}
Refer to the Appendix.
\end{proof}

\begin{figure}[tb] \centering

\subfigure[Point locations.]{\includegraphics[scale=0.34]{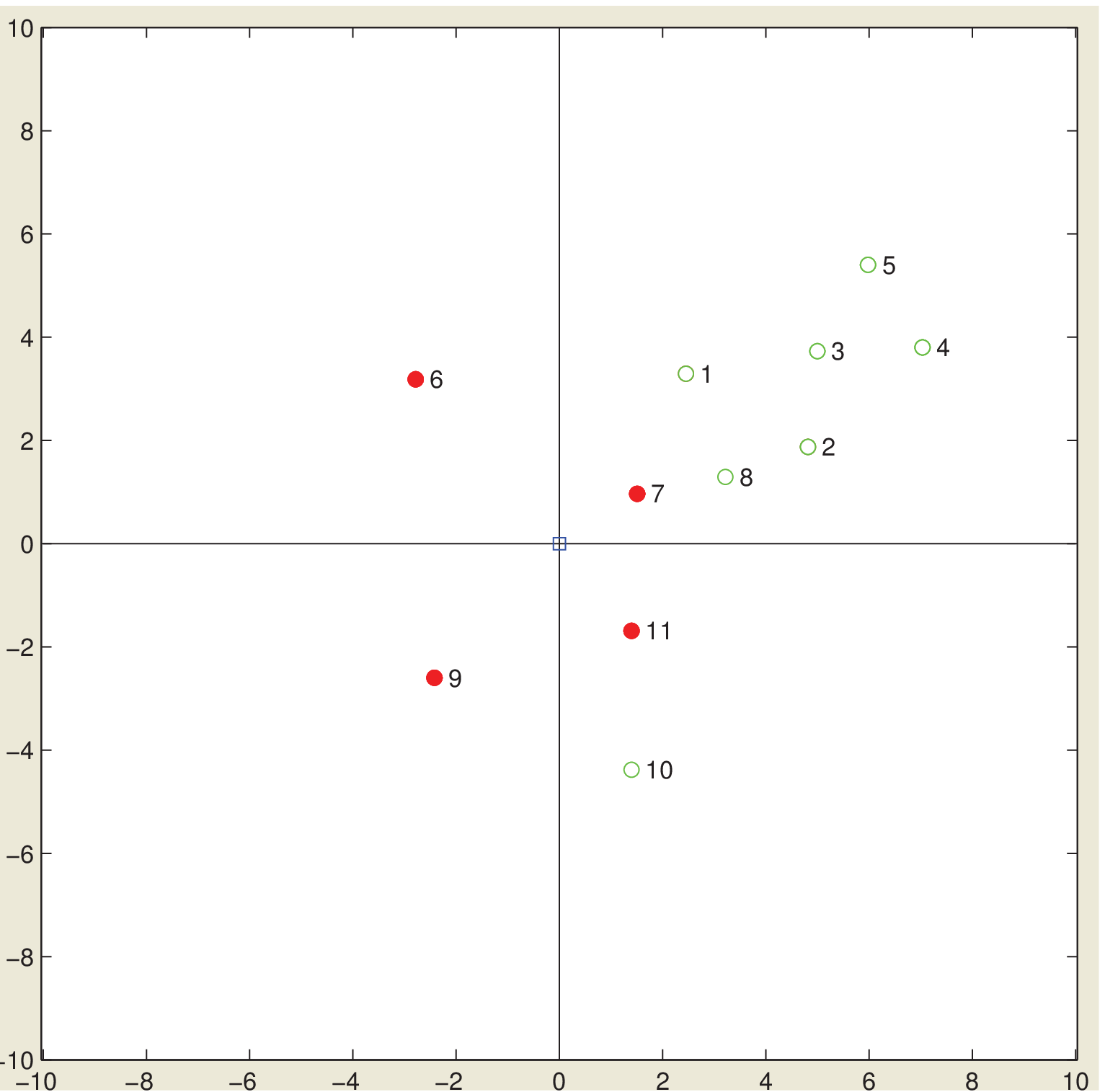}} \hfil \subfigure[Voronoi cell.]{\includegraphics[scale=0.38]{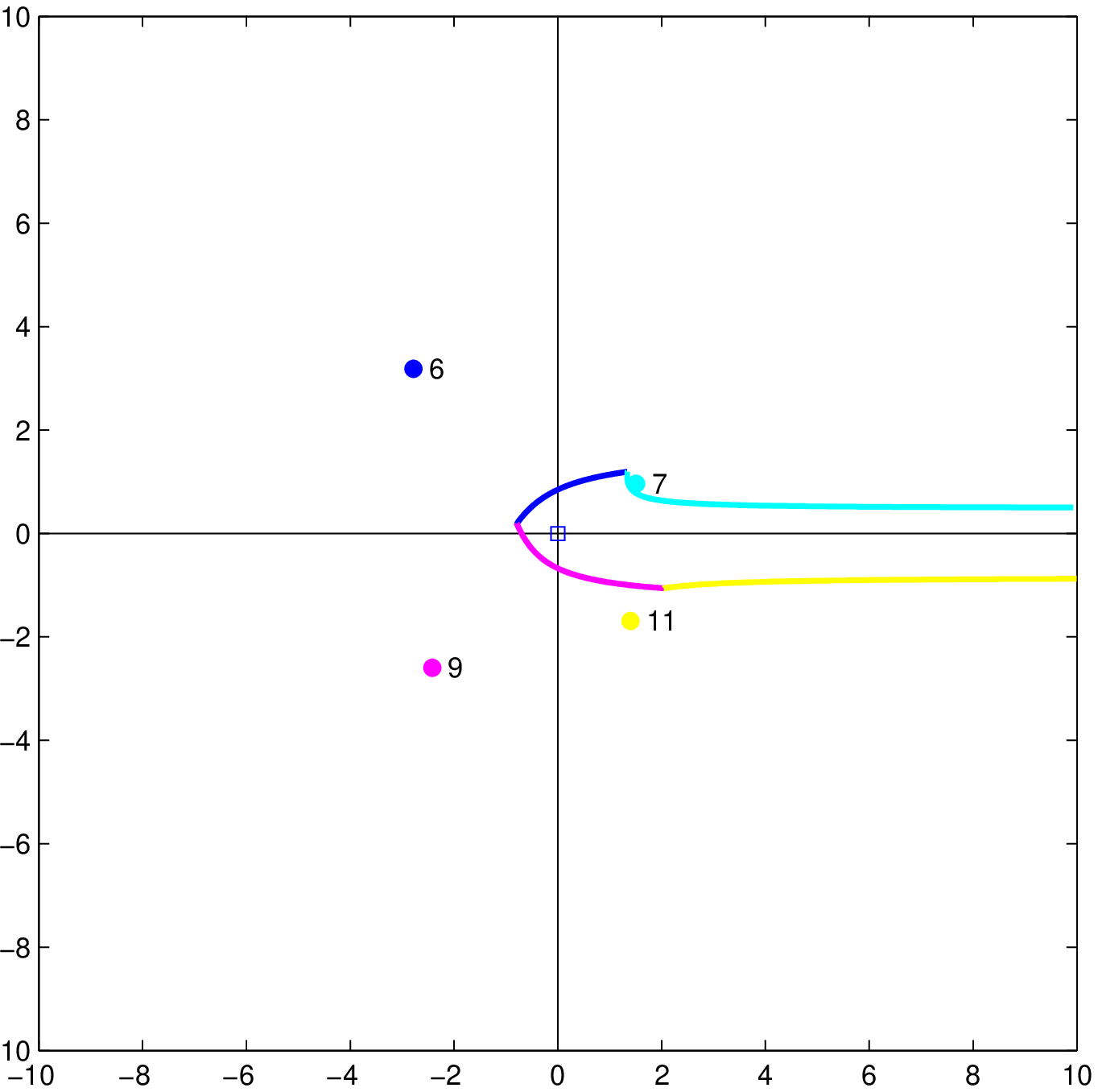}}
\caption{A fixed point $p^1 = (0~0)^T$ and a set $\mathds{P}$ of $11$ points in (a), and the (bounded) Voronoi cell of the point $p^1$ in (b). $\mathcal{N}_G(p^1)$ consists of points with the red color in (a), which are also Voronoi neighbors (for this example). In (b), the same color is used for plotting the point $i$ and the Voronoi boundary between point $i$ and $p^1$.} \label{fig10} \vspace{-15pt}
\end{figure}

\begin{example}
Let $p^1 = (0~0)^T$, and generate a set $\mathds{P}$ of $11$ points that satisfy Assumption~\ref{assumption} in the square $[-10, 10]\times [-10, 10]$ as shown in Fig.~\ref{fig10}(a). It can be verified that $\mathcal{N}_G(p^1)$ consists of points $6, 7, 9, 11$ (which are highlighted using the red color). By calculating the Voronoi cell of the point $p^1$ (the bounded version is shown in Fig.~\ref{fig10}(b)), the Voronoi neighbors of $p^1$ are points $6, 7, 9, 11$, which are the same as the set of points $\mathcal{N}_G(p^1)$ (for this example). \label{example:upper_bound} \hfill $\diamondsuit$
\end{example}

%

\section{Calculation of the Upper Bound: Algorithms} \label{section6}
In this section, we first discuss how to calculate the upper bound (on the set of Voronoi neighbors) in Theorem~\ref{prop:Voronoi_neighbor} when the set of generators is fixed, and then propose algorithms to deal with dynamically changing generators.

\subsection{Algorithms for Calculating the Upper Bound: Static Case} \label{section6_static}
Since the set of Voronoi neighbors is important for solving Problem~\ref{problem:distributed_calculation}, it is necessary to develop algorithms to calculate the upper bound in Theorem~\ref{prop:Voronoi_neighbor}. One simple algorithm is given in Algorithm~\ref{algorithm_simple}. The algorithm just checks the condition in Theorem~\ref{prop:Voronoi_neighbor}. More specifically, the variable $\mathit{sign}$ indicates whether the point $p^i$ belongs to $\mathcal{N}_{G}(p^1)$: $1$ if it does and $0$ otherwise. Steps 4-6 verify if $p^i$ is dominated by some point $p^j$: if it is, then the algorithm exits the inner \textbf{for} loop and does not add $p^i$ into $\mathcal{N}_{G}(p^1)$; otherwise (i.e., $\mathit{sign}$ is never set to be $0$), the algorithm adds $p^i$ to $\mathcal{N}_{G}(p^1)$. It can be verified that the algorithm has complexity $\mathcal{O}(n^2)$. The complexity $\mathcal{O}(n^2)$ is tight for the scenario in which none of $p^2, p^3, ..., p^n$ dominates any other point, i.e., $\mathcal{N}_{G}(p^1) = \mathds{P}$.

\begin{algorithm}[tb] \small
\caption{Simple Upper Bound Calculation} \label{algorithm_simple}

\begin{algorithmic}

\REQUIRE A fixed point $p^1$ and a set of points $\mathds{P} = \{p^2, p^3, ..., p^n\}$ satisfying Assumption~\ref{assumption}

\ENSURE $\mathcal{N}_{G}(p^1)$ \vspace{2pt}

\end{algorithmic}

\begin{algorithmic}[1]

\STATE Initialize $\mathcal{N}_{G}(p^1) = \emptyset$;

\FOR{$i = 2, 3, \dots n$}

\STATE Let $\mathit{sign} = 1$;

\FOR{$j = 2, 3, \dots n$}

\STATE If $j \neq i$ and $p^j \succ p^i$, set $\mathit{sign} = 0$ and exit the inner \textbf{for} loop;

\ENDFOR

\STATE If $\mathit{sign} = 1$, set $\mathcal{N}_{G}(p^1) = \mathcal{N}_{G}(p^1) \cup \{p^i\}$;

\ENDFOR

\STATE Output $\mathcal{N}_{G}(p^1)$.
\end{algorithmic}
\end{algorithm}

One natural question to ask is whether there exists more efficient algorithms to calculate the upper bound. Note that the conditions in Theorem~\ref{case_1} require that $y_{p^3} \geq y_{p^2}$. Therefore, if we first sort the points according to $y$ coordinates, potentially we can obtain a faster algorithm.

Given the set of points $\mathds{P} = \{p^2, p^3, ..., p^n\}$, we first divide the set of points into three groups: points of which the $y$ coordinate is larger than $y_{p^1}$ (denoted as $\mathds{P}_+$), points of which the $y$ coordinate is the same as $y_{p^1}$ (denoted as $\mathds{P}_0$), and points of which the $y$ coordinate is smaller than $y_{p^1}$ (denoted as $\mathds{P}_-$).

We first focus on points in $\mathds{P}_+$. If $\mathds{P}_+$ is not empty, we sort the points in $\mathds{P}_+$ according to the ascending order of $y$ coordinates and according to the ascending order of $x$ coordinates for points that have the same $y$ coordinate. Suppose there are $m$ points in $\mathds{P}_+$ and the sorted point sequence is $p_1, p_2, ..., p_m$. We use the point $\mathit{anchor}$ to track the angle formed by the positive $x$ axis and the ray from $p^1$ to the point (if the point is $p$, we use $\angle p$ to denote this angle). Since $p_1$ has the smallest $y$, it must belong to $\mathcal{N}_G(p^1)$ due to Theorem~\ref{case_1}, and $\mathit{anchor}$ is initialized as $p_1$. Now we consider $p_2$: if $p_2$ is not dominated by the point $\mathit{anchor}$, then add $p_2$ into $\mathcal{N}_G(p^1)$ and update $\mathit{anchor}$ with $p_2$; otherwise, do nothing. We repeat this procedure for points $p_3, p_4, ..., p_m$.

For points in $\mathds{P}_-$, we can simply change the sign of the $y$ coordinates so that we can use the procedure for $\mathds{P}_+$ due to Corollary~\ref{case_2}. For points in $\mathds{P}_0$, we add the point with the smallest $x$ coordinate into $\mathcal{N}_G(p^1)$. The detailed algorithm is given in Algorithm~\ref{algorithm_efficient}.

The correctness of the algorithm can be proved as below. The reason why $\mathds{P}$ can be divided into three subsets $\mathds{P}_+$, $\mathds{P}_0$, and $\mathds{P}_-$ is that for any point $p$ in each subset, it can only potentially be dominated by points in that subset due to the dominance characterizations in Theorem~\ref{case_1}, Corollary~\ref{case_2}, and part (b) of Proposition~\ref{prop:special_case}. For points in $\mathds{P}_+$, we obtain $p_1, p_2, ..., p_m$ after sorting.  Due to the condition in Theorem~\ref{case_1}, a point $p_i$ can only be dominated by points that have indices smaller than $i$. Thus, $p_1$ must belong to $\mathcal{N}_G(p^1)$ because it has the smallest index, and we initialize $\mathit{anchor}$ with $p_1$. When considering the point $p_i$ with $i = 2, 3, ..., m$, the variable $\mathit{anchor}$ keeps track of the point $p$ that has formed the largest angle $\angle p$ so far. If $p_i$ is dominated by the point $\mathit{anchor}$, it cannot belong to $\mathcal{N}_{G}(p^1)$. If $p_i$ is not dominated by the point $\mathit{anchor}$, then $p_i$ cannot be dominated by any point in $\mathds{P}_+$ (and we add $p_i$ into $\mathcal{N}_{G}(p^1)$ and  update $\mathit{anchor}$ with $p_i$); this is because

\begin{itemize}

\item[i)] $p_i$ can only be potentially dominated by points that have indices smaller than $i$,

\item[ii)] for any point $p$ with its index smaller than $i$, we have $y_{p_i} \geq y_p$ so we only need to check the second condition in Theorem~\ref{case_1}; however, the second condition essentially just compares the angle $\angle p$ with the angle $\angle p_i$. If there exists a point $p_k$ such that $\angle p_k \geq \angle p_i$, then $p_i$ is dominated by $p_k$. Since the variable $\mathit{anchor}$ keeps track of the point which forms the largest angle so far, it is equivalent to compare the angle $\angle \mathit{anchor}$ with $\angle p_i$, which is the same as checking if $p_i$ is dominated by the point $\mathit{anchor}$.
\end{itemize}
Because the condition in Corollary~\ref{case_2} is symmetric to the condition in Theorem~\ref{case_1}, we can just change the sign of the $y$ coordinates for points in $\mathds{P}_-$ and apply the procedure for points in $\mathds{P}_+$. For points in $\mathds{P}_0$, there is only one point (namely, the one with the smallest $x$ coordinate) that cannot be dominated by any other point due to part (b) of Proposition~\ref{prop:special_case}.

It can be verified that the algorithm has complexity $\mathcal{O}(n \log n)$ because the best sorting algorithms (e.g., merge sort) have complexity $\mathcal{O}(n \log n)$ and finding points in $\mathcal{N}_{G}(p^1)$ using the sorted point sequence has complexity $\mathcal{O}(n)$.

\begin{algorithm}[tb] \small
\caption{Efficient Upper Bound Calculation} \label{algorithm_efficient}

\begin{algorithmic}

\REQUIRE A fixed point $p^1$ and a set of points $\mathds{P} = \{p^2, p^3, ..., p^n\}$ satisfying Assumption~\ref{assumption}

\ENSURE $\mathcal{N}_{G}(p^1)$ \vspace{2pt}

\end{algorithmic}

\begin{algorithmic}[1]

\STATE Initialize $\mathcal{N}_{G}(p^1) = \emptyset$;

\STATE Divide the set of points in $\mathds{P}$ into three subsets $\mathds{P}_+$, $\mathds{P}_0$, and $\mathds{P}_-$;

\IF{$|\mathds{P}_+| > 0$}

\STATE Sort the points in $\mathds{P}_+$ into $p_1, p_2, ..., p_m$ (with $m = |\mathds{P}_+|$) according to the ascending order of $y$ coordinates and according to the ascending order of $x$ coordinates for points that have the same $y$ coordinate;

\STATE Add $p_1$ into $\mathcal{N}_{G}(p^1)$ and set $\mathit{anchor} = p_1$;

\FOR{$i = 2, 3, \dots m$}

\STATE If $p_i$ is not dominated by the point $\mathit{anchor}$, add $p_i$ into $\mathcal{N}_{G}(p^1)$ and set $\mathit{anchor} = p_i$;

\ENDFOR

\ENDIF

\STATE If $|\mathds{P}_0| > 0$, add the point with the smallest $x$ coordinate in $\mathds{P}_0$ into $\mathcal{N}_{G}(p^1)$;

\STATE If $|\mathds{P}_-| > 0$, let $\mathds{P}_-' = \{p~|~(x_p~-y_p)^T \in \mathds{P}_- \}$, and apply the procedure in Steps~3-9 by replacing $\mathds{P}_+$ with $\mathds{P}_-'$ and adding $(x_{p_i}~-y_{p_i})^T$ into $\mathcal{N}_G(p^1)$;

\STATE Output $\mathcal{N}_{G}(p^1)$.
\end{algorithmic}
\end{algorithm}

\subsection{Algorithms for Calculating the Upper Bound: Dynamic Case}
In this subsection, we study how to calculate the upper bound if the set of generators $\mathds{P}$ constantly changes. If we apply the algorithms in Section~\ref{section6_static} to changing generators, the calculation has to be done for any single change in the generator locations. Here, we introduce a dominance graph, which builds up on the dominance relation and can be used to calculate the upper bound more efficiently. We first look at another property of the dominance relation.

\begin{proposition} (Transitivity of Dominance)
Given a fixed point $p^1$, and three points $p^2, p^3, p^4$ that are different from $p^1$, if $p^2 \succ p^3$ and $p^3 \succ p^4$, then $p^2 \succ p^4$. \label{prop:transitivity}
\end{proposition}

\begin{proof}
Refer to the Appendix.
\end{proof}

Under Assumption~\ref{assumption}, the dominance relation is a partial order since it is
\begin{itemize}
\item Reflexive because $\forall p \in \mathds{P}$, $p \succ p$;
\item Antisymmetric because $p^i \succ p^j$ and $p^j \succ p^i$ imply that $p^i = p^j$ as shown in Proposition~\ref{prop:antisymmetry};
\item Transitive because $p^i \succ p^j$ and $p^j \succ p^k$ imply that $p^i \succ p^k$ as shown in Proposition~\ref{prop:transitivity}.
\end{itemize}
Based on this dominance partial order, the set of points $\mathds{P}$ induces a directed graph, which we call a dominance graph.

\begin{definition}
Given a fixed point $p^1$ and a set of points $\mathds{P} = \{p^2, p^3, ..., p^n\}$ satisfying Assumption~\ref{assumption}, we define the dominance graph induced by $\mathds{P}$ as $\mathds{G}(p^1, \mathds{P}) = (V~E)$, where i) $V = \mathds{P}$ is the set of vertices, and ii) $E$ is the set of directed edges and there is a directed edge $e \in E$ from $p^i \in V$ to $p^j \in V$ if $p^i \succ p^j$. In addition, $p^i$ is called a parent of $p^j$, and $p^j$ is called a child of $p^i$.
\end{definition}

Since a point can dominate (or be dominated by) multiple points in $\mathds{P}$, there could be multiple output (or input) edges for this point in the dominance graph. However, there is no cycle in dominance graphs as shown in the following proposition, which can be proved via contradiction.

\begin{proposition}
Given a fixed point $p^1$ and a set of points $\mathds{P} = \{p^2, p^3, ..., p^n\}$ satisfying Assumption~\ref{assumption}, the dominance graph $\mathds{G}(p^1, \mathds{P})$ is acyclic. \label{prop:acyclic}
\end{proposition}

\begin{proof}
Refer to the Appendix.
\end{proof}

Therefore, dominance graphs are directed acyclic graphs. Given a dominance graph, we are interested in finding points that are not dominated by any other point.

\begin{definition}
Given a fixed point $p^1$ and a set of points $\mathds{P} = \{p^2, p^3, ..., p^n\}$ satisfying Assumption~\ref{assumption}, a point $p \in \mathds{P}$ is a neighbor of $p^1$ in the dominance graph $\mathds{G}(p^1, \mathds{P})$ if there does not exist another point $p' \in \mathds{P}$ that is different from $p$ and dominates $p$. \label{def:neighbor_d_graph}
\end{definition}

Based on Definition~\ref{def:NGp1}, $\mathcal{N}_{G}(p^1)$ is exactly the set of neighbors of $p^1$ in the dominance graph $\mathds{G}(p^1, \mathds{P})$. The importance of neighbors in the dominance graph is that only neighbors of $p^1$ matter when $p^1$ calculates its own Voronoi cell given the set of points $\mathds{P}$, as shown in Theorem~\ref{prop:Voronoi_neighbor}. It can be verified that the following result holds.
\begin{proposition}
Given a fixed point $p^1$ and a set of points $\mathds{P} = \{p^2, p^3, ..., p^n\}$ satisfying Assumption~\ref{assumption}, $p$ is a neighbor of $p^1$ in the dominance graph iff the in-degree of $p$ in the dominance graph $\mathds{G}(p^1, \mathds{P})$ is $0$. \label{prop:neighbor}
\end{proposition}

\begin{figure*}[tb] \centering
\subfigure[Dominance graph $\mathds{G}(p^1, \mathds{P})$ for the setting in Fig.~\ref{fig10}(a).]{\includegraphics[scale=0.32]{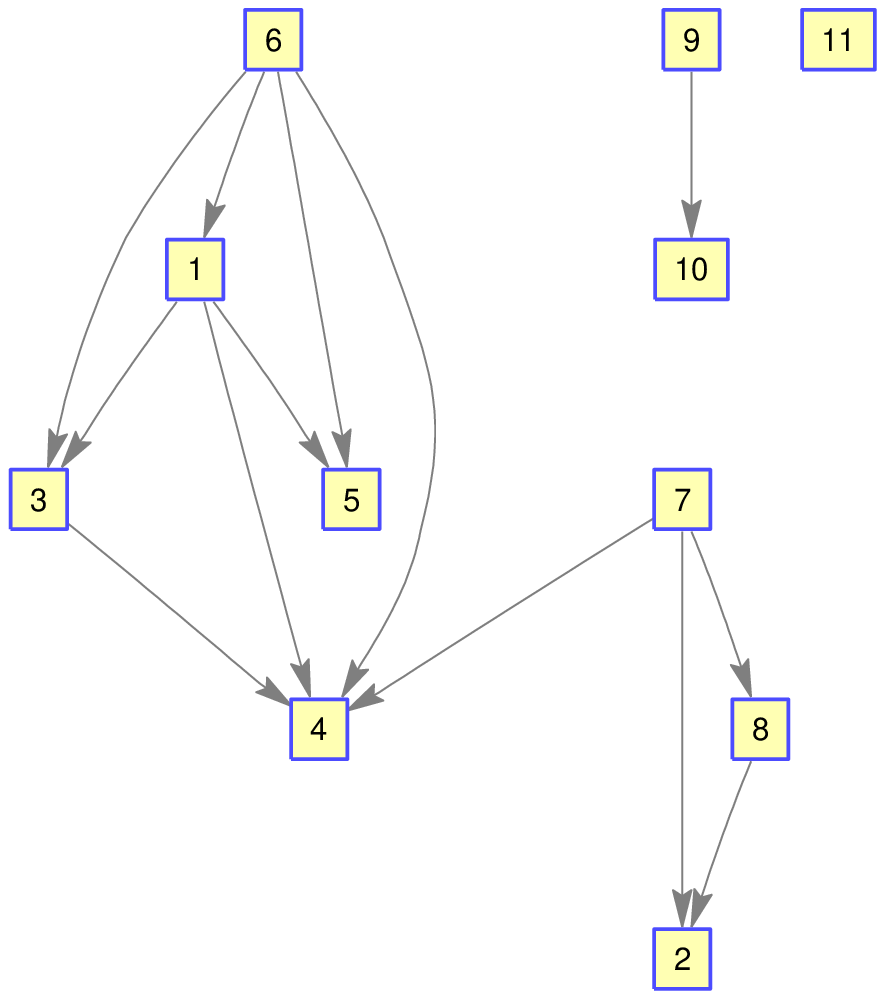}} \hfil
\subfigure[Inserting points $12$ and $13$ in Fig.~\ref{fig10}(a).]{\includegraphics[scale=0.28]{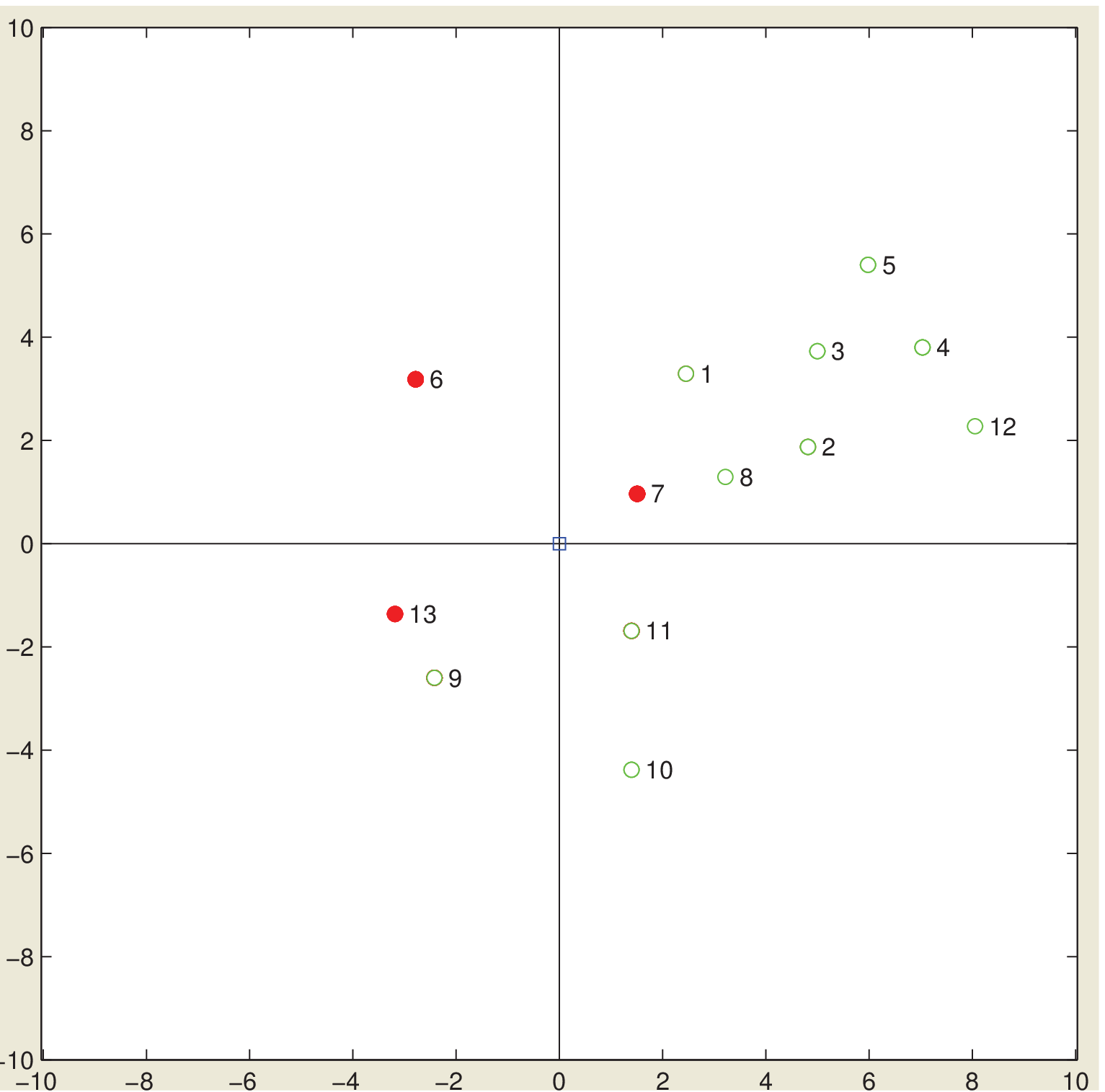}} \hfil
\subfigure[Dominance graph after inserting point $12$, which does not become a neighbor.]{\includegraphics[scale=0.32]{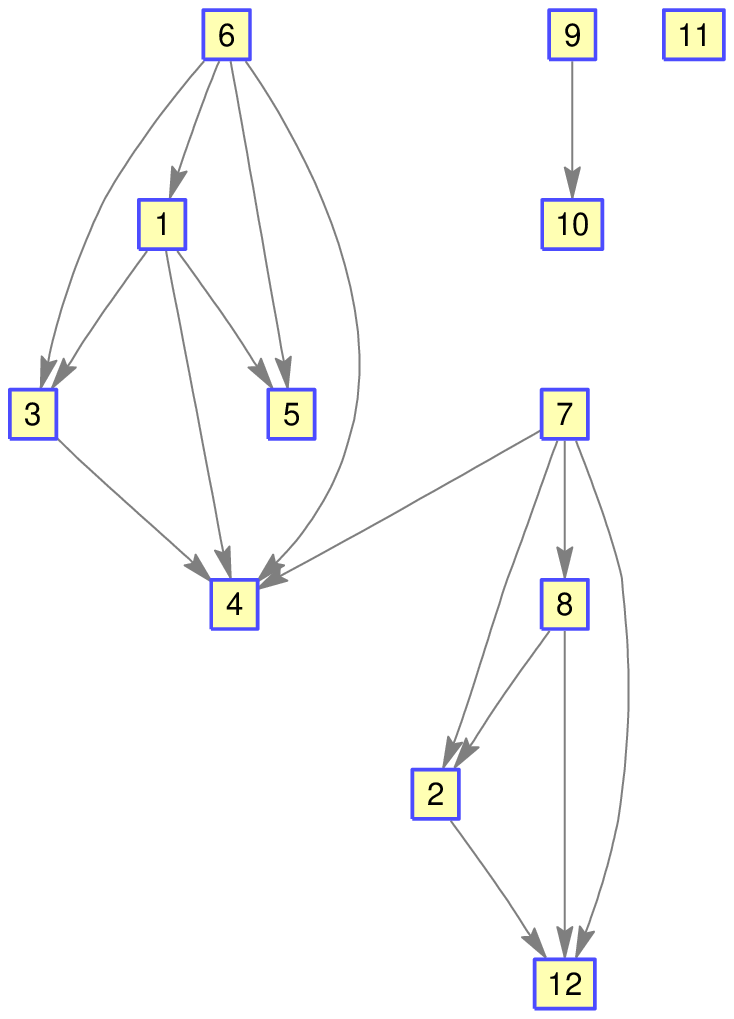}} \hfil \subfigure[Dominance graph after inserting point $13$, which becomes a neighbor.]{\includegraphics[scale=0.32]{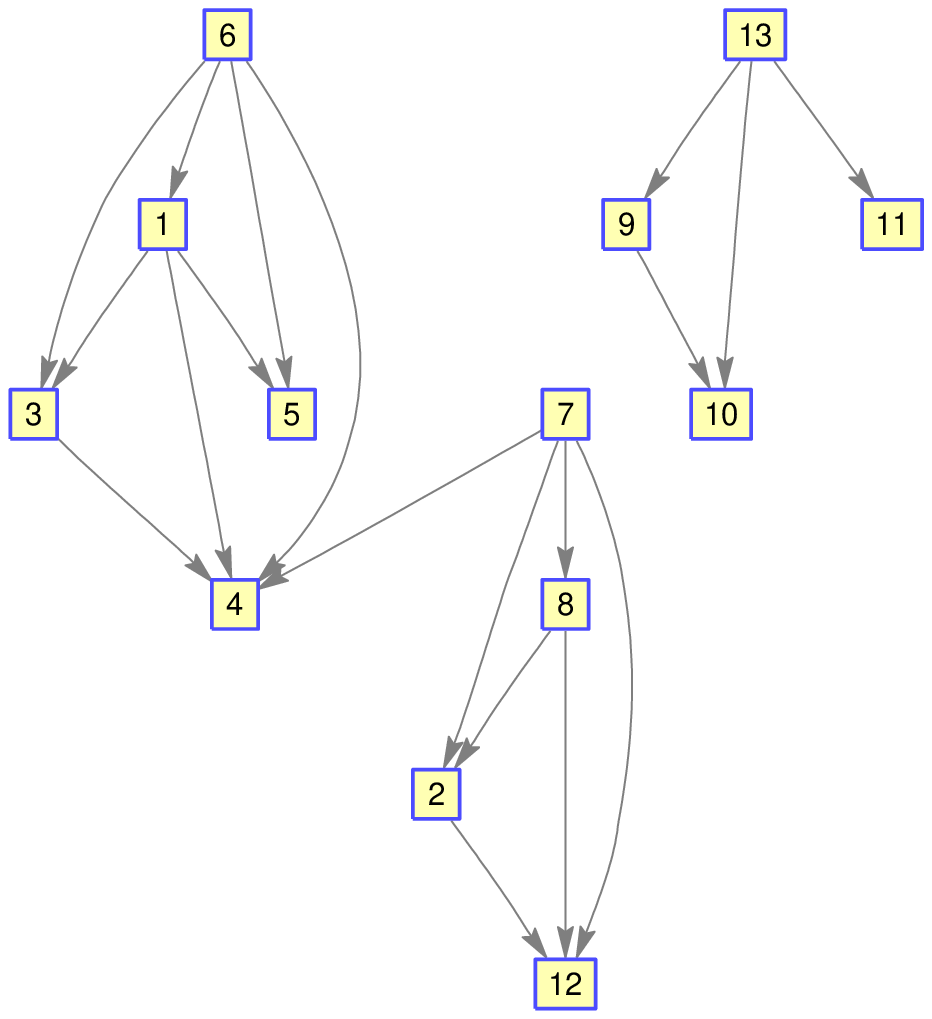}} \hfil
\subfigure[Dominance graph after deleting point $5$, which is not a neighbor.]{\includegraphics[scale=0.32]{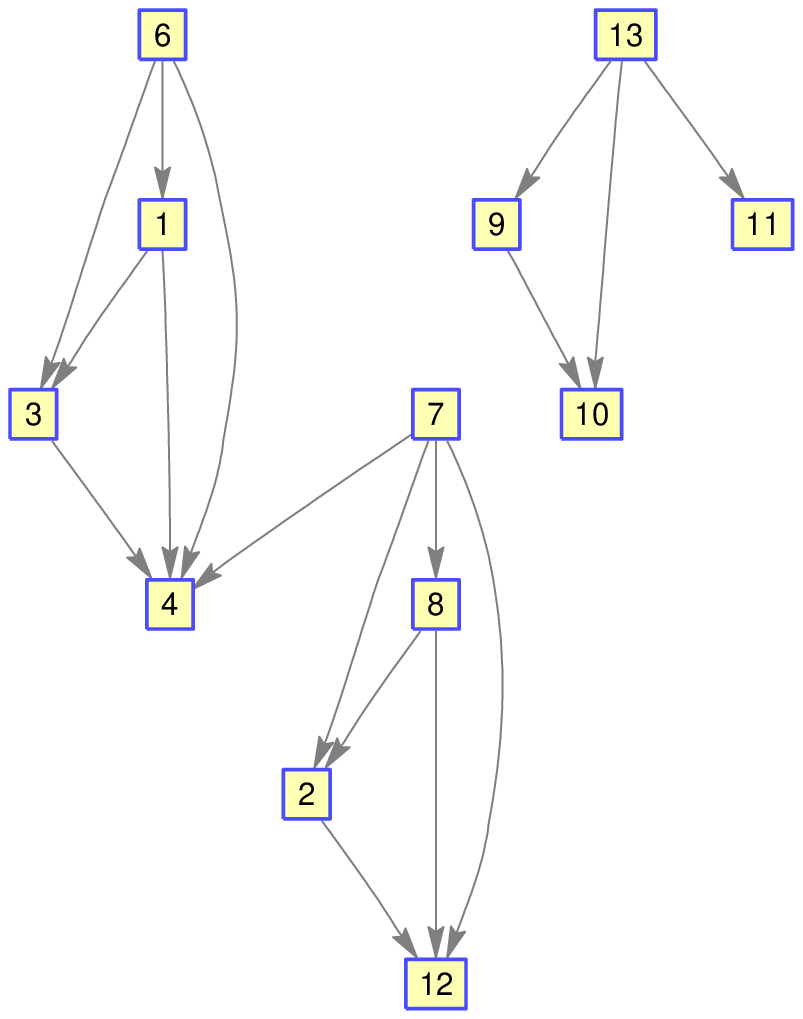}}
\hfil
\subfigure[Dominance graph after deleting point $7$, which is a neighbor.]{\includegraphics[scale=0.30]{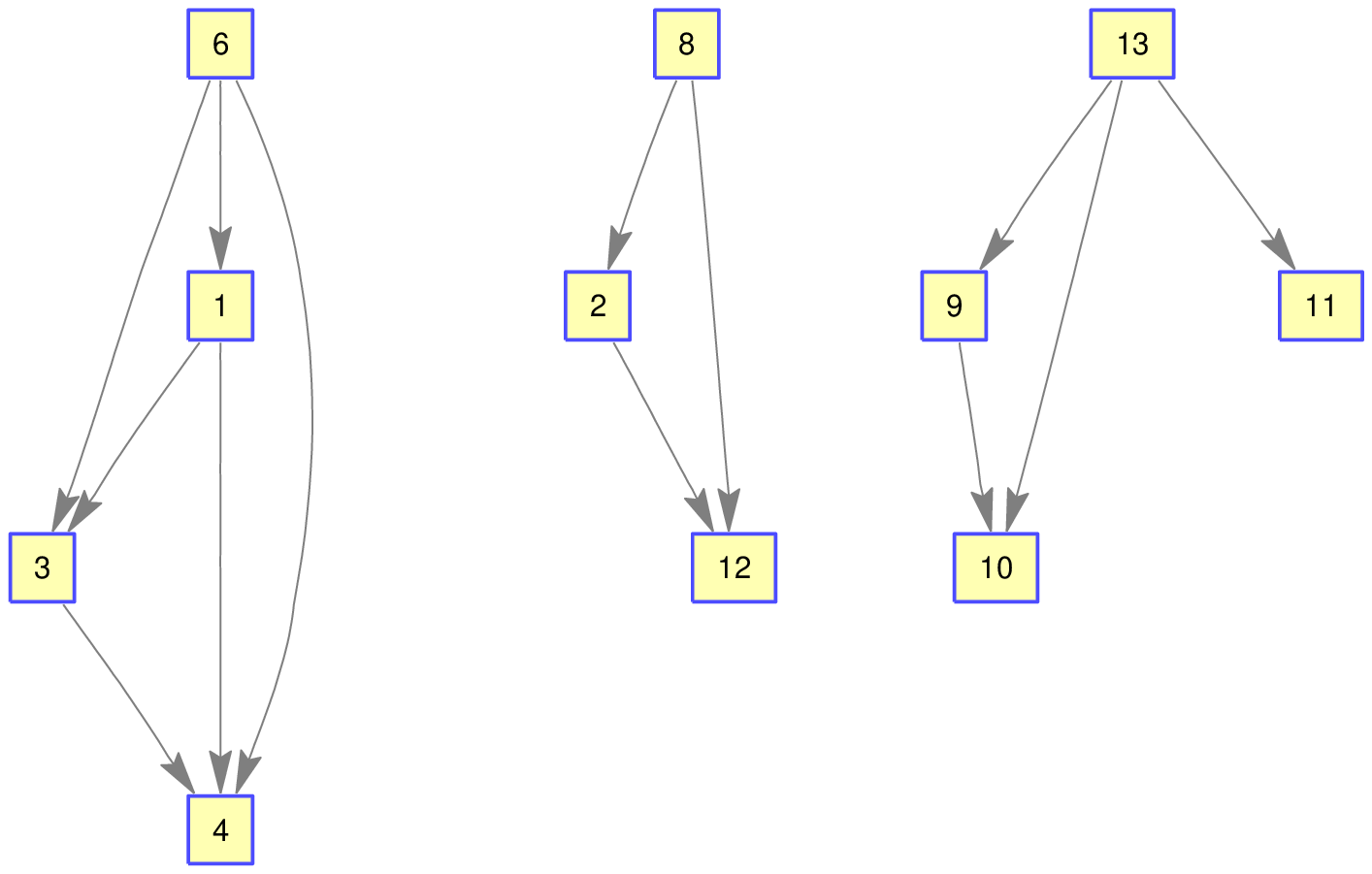}}
\caption{Dynamic maintenance of dominance graph.} \label{fig12} \vspace{-15pt}
\end{figure*}

\begin{example} We still consider the setting in Example~\ref{example:upper_bound}: $p^1 = (0~0)^T$, and the set of $11$ points $\mathds{P}$ is shown in Fig.~\ref{fig10}(a). The corresponding dominance graph is shown in Fig.~\ref{fig12}(a). Based on Proposition~\ref{prop:neighbor}, points $6, 7, 9, 11$ are the neighbors of $p^1$ in the dominance graph. Note that there is no cycle in the dominance graph, which is consistent with Proposition~\ref{prop:acyclic}. \hfill $\diamondsuit$ \label{example:dominance_graph}
\end{example}



Now we study how to dynamically maintain the dominance graph when inserting or deleting points. We assume that there is an upper bound $K$ on the number of points (naturally, the number of mobile vehicles serves as this upper bound $K$). For each point $p^i \in \mathds{P}$, we have the vertex $V_{p^i}$ in the dominance graph, and use the following fields to keep track of the vertex: i) $\textrm{ID}$: $i$ with $2 \leq i \leq K$; ii) $x$: the $x$ coordinate of $p^i$; iii) $y$: the $y$ coordinate of $p^i$; iv) \textit{Parent}: a data array of dimension $K$. \textit{Parent}$(k) = 1$ if there is an edge from $V_{p^k}$ to $V_{p^i}$, $0$ otherwise; v) \textit{Child}: a data array of dimension $K$. \textit{Child}$(k) = 1$ if there is an edge from $V_{p^i}$ to $V_{p^k}$, $0$ otherwise; vi) \textit{No\_of\_parent}: the number of parents of vertex $V_{p^i}$. In addition, there is a data array \textit{List\_of\_neighbor} of dimension $K$ for $p^1$ to keep track of its neighbors, i.e., \textit{List\_of\_neighbor}$(k) = 1$ if $p^k$ is a neighbor of $p^1$ and $0$ otherwise.

When inserting a point, the point can affect the child and parent fields of the existing vertices, can make a neighbor invalid, and itself can become a new neighbor. The details are given in Algorithm~\ref{algorithm_insertion}. The input is a dominance graph $\mathds{G}(p^1, \mathds{P})$ with vertices $V_{p^2}, ..., V_{p^n}$, a list of neighbors, and a new point $p^{n+1}$, and the output is the dominance graph $\mathds{G}(p^1, \mathds{P} \cup \{p^{n+1}\})$ and the updated list of neighbors. Step~1 initializes the vertex $V_{p^{n+1}}$. Steps~2-11 update the fields of vertices $V_{p^2}, V_{p^3}, ..., V_{p^{n+1}}$. At Step~3, $p^{n+1}$ is checked against $p^i$, and there are three cases:
 \begin{itemize}
 \item $p^i$ dominates $p^{n+1}$. Then $p^i$ is a parent of $p^{n+1}$, and the number of parents of $p^{n+1}$ is increased by $1$. This corresponds to Step~4;
 \item $p^{n+1}$ dominates $p^i$. Then $p^i$ is a child of $p^{n+1}$, the number of parents of $p^i$ is increased by $1$, and $p^i$ cannot be a neighbor of $p^1$. This corresponds to Step~5.
 \item Neither $p^i$ nor $p^{n+1}$ dominates. In this case, no change is necessary.
 \end{itemize}
 Steps~7-11 examine whether $p^{n+1}$ itself is a neighbor of $p^1$. It can be verified that the insertion algorithm has complexity $\mathcal{O}(n)$.  Algorithm~\ref{algorithm_insertion} can start with $n = 1$ (or $n \geq 2$), i.e., when the dominance graph is empty.


\begin{algorithm}[tb] \small
\caption{Dynamic Point Insertion} \label{algorithm_insertion}

\begin{algorithmic}

\REQUIRE Dominance graph $\mathds{G}(p^1, \mathds{P})$ with each vertex $V_p$ for $p \in P$ being represented using the fields \textit{ID}, $x$, $y$, \textit{Parent}, \textit{Child}, and \textit{No\_of\_parent}, a global data structure \textit{List\_of\_neighbor}, and a new point $p^{n+1}$ such that Assumption~\ref{assumption} holds for the set of points $\mathds{P} \cup \{p^{n+1}\}$

\ENSURE Dominance graph $\mathds{G}(p^1, \mathds{P} \cup \{p^{n+1}\})$ and updated \textit{List\_of\_neighbor} \vspace{2pt}

\end{algorithmic}

\begin{algorithmic}[1]

\STATE Initialize $V_{p^{n+1}}$ with \textit{ID} being $n+1$, $x$ being $x_{p^{n+1}}$, $y$ being $y_{p^{n+1}}$, \textit{Parent} and \textit{Child} being vectors of all zeros, and \textit{No\_of\_parent} being $0$;

\FOR{$i = 2, 3, \dots n$}

\STATE Check if $p^i$ dominates $p^{n+1}$ or $p^{n+1}$ dominates $p^i$ using the results in Theorem~\ref{case_1}, Corollary~\ref{case_2}, and Proposition~\ref{prop:special_case};

\STATE If $p^i$ dominates $p^{n+1}$, increase $V_{p^{n+1}}.$\textit{No\_of\_parent} by $1$, and set $V_{p^i}.$\textit{Child}$(n+1) = 1$ and $V_{p^{n+1}}.$\textit{Parent}$(i) = 1$;

\STATE If $p^{n+1}$ dominates $p^{i}$, increase $V_{p^i}.$\textit{No\_of\_parent} by $1$, and set $V_{p^i}.$\textit{Parent}$(n+1) = 1$, \textit{List\_of\_neighbor}$(i) = 0$, and $V_{p^{n+1}}.$\textit{Child}$(i) = 1$;

\ENDFOR

\IF{$V_{p^{n+1}}.$\textit{No\_of\_parent}$= 0$}

\STATE \textit{List\_of\_neighbor}$(n+1) = 1$;

\ELSE

\STATE \textit{List\_of\_neighbor}$(n+1) = 0$.

\ENDIF

\end{algorithmic}
\end{algorithm}

When deleting a point, the point can affect the child and parent fields of other vertices, and can create new neighbors. The details are given in Algorithm~\ref{algorithm_deletion}. The input is a dominance graph $\mathds{G}(p^1, \mathds{P})$ with vertices $V_{p^2}, ..., V_{p^n}$, a list of neighbors, and a point $p^{j}$ to delete, and the output is the dominance graph $\mathds{G}(p^1, \mathds{P} \setminus \{p^{j}\})$ and the updated list of neighbors. Steps~1-5 update the fields of vertices $V_{p^2}, V_{p^3}, ..., V_{p^{j-1}}, V_{p^{j+1}}, ..., V_{p^n}$. If $p^i$ is a parent of $p^j$, then remove $p^j$ from $p^i$'s child list; this is done in Step~2. If $p^i$ is a child of $p^j$, then remove $p^j$ from $p^i$'s parent list, decrease the number of parents of $p^i$ by $1$: if the number of parents of $p^i$ is $0$, set $p^i$ to be a neighbor of $p^1$; this is done in Step~4. If $p^j$ is a neighbor of $p^1$, Step~6 removes $p^j$ from the list of neighbors of $p^1$. It can be verified that the deletion algorithm has complexity $\mathcal{O}(n)$. Algorithm~\ref{algorithm_deletion} can start with $n=2$ (or $n \geq 3$), i.e., when there is only one vertex in the dominance graph.

\begin{algorithm}[tb] \small
\caption{Dynamic Point Deletion} \label{algorithm_deletion}

\begin{algorithmic}

\REQUIRE Dominance graph $\mathds{G}(p^1, \mathds{P})$ with each vertex $V_p$ for $p \in P$ being represented using the fields \textit{ID}, $x$, $y$, \textit{Parent}, \textit{Child}, and \textit{No\_of\_parent}, a global data structure \textit{List\_of\_neighbor}, and a point $p^{j} \in \mathds{P}$ to delete, where $j \in \{2, 3, \dots, n\}$

\ENSURE Dominance graph $\mathds{G}(p^1, \mathds{P} \setminus \{p^{j}\})$ and updated \textit{List\_of\_neighbor} \vspace{2pt}

\end{algorithmic}

\begin{algorithmic}[1]

\FOR{$i = 2, 3, \dots n$}

\STATE If $V_{p^j}.$\textit{Parent}$(i) = 1$, set $V_{p^i}.$\textit{Child}$(j) = 0$;

\STATE If $V_{p^j}.$\textit{Child}$(i) = 1$, set $V_{p^i}.$\textit{Parent}$(j) = 0$ and decrease $V_{p^i}.$\textit{No\_of\_parent} by $1$;

\STATE If $V_{p^i}.$\textit{No\_of\_parent}$=0$, set \textit{List\_of\_neighbor}$(i) = 1$;

\ENDFOR

\STATE Set \textit{List\_of\_neighbor}$(j) = 0$, and delete the vertex $V_{p^j}$.

\end{algorithmic}
\end{algorithm}

To obtain the set of neighbors of $p^1$ in the dominance graph at any time, we just need to check the nonzero entries in \textit{List\_of\_neighbor}.

\begin{example}
We still use the point $p^1$ and the set of points $\mathds{P}$ in Example~\ref{example:upper_bound}. Now we consider inserting points. We first insert point $12$ (the relative position of which is shown in Fig.~\ref{fig12}(b)). The dominance graph in Fig.~\ref{fig12}(a) becomes the graph in Fig.~\ref{fig12}(c); here, point $12$ is not a neighbor of $p^1$ in its dominance graph (since it is dominated by points $2$, $7$ and $8$), which implies that the Voronoi cell of $p^1$ remains the same following Theorem~\ref{prop:Voronoi_neighbor}. Next we insert point $13$ (the relative position of which is shown in Fig.~\ref{fig12}(b)). The dominance graph in Fig.~\ref{fig12}(c) becomes the graph in Fig.~\ref{fig12}(d); here, point $13$ becomes a neighbor of $p^1$ in its dominance graph. In this case, points $6, 7, 13$ are neighbors of $p^1$ in its dominance graph, and can also be verified to be Voronoi neighbors. Now we consider deleting points. We first delete point $5$ (the relative position of which is shown in Fig.~\ref{fig12}(b)). The dominance graph in Fig.~\ref{fig12}(d) becomes the graph in Fig.~\ref{fig12}(e). Since the set of neighbors in dominance graph does not change after deleting point $5$, the Voronoi cell of $p^1$ remains the same as the one after inserting point $13$ following Theorem~\ref{prop:Voronoi_neighbor}. Next we delete point $7$ (the relative position of which is shown in Fig.~\ref{fig12}(b)). The dominance graph in Fig.~\ref{fig12}(e) becomes the graph in Fig.~\ref{fig12}(f). Since point $7$ is a neighbor of $p^1$ in its dominance graph, point $8$ becomes a neighbor of $p^1$ in its dominance graph after deleting point $7$, as shown in Fig.~\ref{fig12}(f). In this case, points $6, 8, 13$ are neighbors of $p^1$ in its dominance graph, and can also be verified to be Voronoi neighbors. Therefore, when a point is inserted or deleted, the Voronoi cell could stay the same provided that the set of neighbors in the dominance graph remains the same. In contrast, we have to consider all points again after insertion or deletion if without the dynamically maintained dominance graph. This becomes particularly useful when the points are the positions of mobile vehicles. \hfill $\diamondsuit$ \label{example:dynamic_maintain}
\end{example}


\section{Simulations} \label{section7}
In this section, we first run simulations to study the number of neighbors in dominance graphs (namely, the cardinality of $\mathcal{N}_{G}(p^1)$), and then propose methods to solve Problem~\ref{problem:distributed_calculation}.

\subsection{Simulations for the Upper Bound}
Since $\mathcal{N}_{G}(p^1)$ is an upper bound on the set of Voronoi neighbors as shown in Theorem~\ref{prop:Voronoi_neighbor}, it is natural to ask the question of how many points there could be in $\mathcal{N}_{G}(p^1)$. While the exact number depends on the specific relative positions of the set of points $\mathds{P}$, we can study the average number of points in $\mathcal{N}_{G}(p^1)$ if the points are generated randomly.

Let us first fix $p^1$ to be $(0~0)^T$, select the number of points $n$ to generate, and then generate each point with uniform distribution in the square $[-1, 1] \times [-1, 1]$ independently from all other points while satisfying Assumption~\ref{assumption}. In other words, the positions of these $n$ points are i.i.d (independent and identically distributed). Given the set of generated $n$ points $\{p^2, ..., p^n, p^{n+1}\}$, we can run either Algorithm~\ref{algorithm_simple} or Algorithm~\ref{algorithm_efficient} to calculate $\mathcal{N}_{G}(p^1)$. Then the quantity that we are interested in is the expected number of points in $\mathcal{N}_{G}(p^1)$.

\begin{figure}[tb] \centering
\includegraphics[scale=0.25]{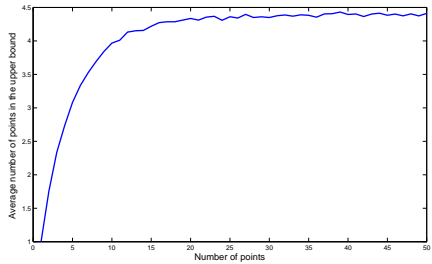}
\caption{Average number of points in $\mathcal{N}_{G}(p^1)$ for randomly generated points over $5000$ trials as a function of the number of points.} \label{fig17}
\end{figure}

\begin{figure*}[tb] \centering
\includegraphics[scale=0.45]{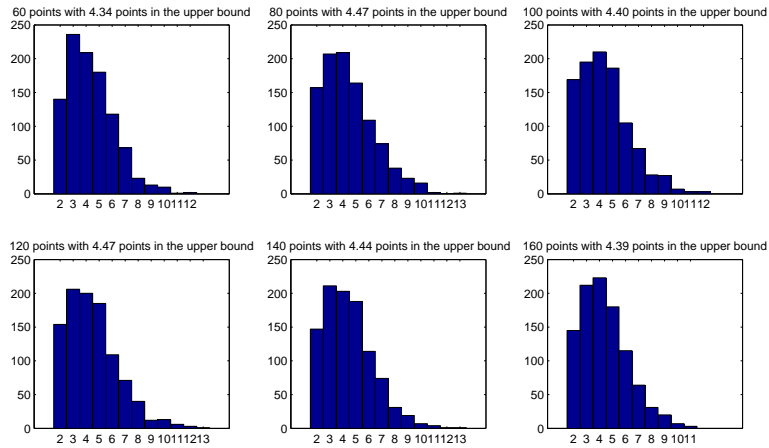}
\caption{Histograms of the number of points (and the average number of points) in $\mathcal{N}_{G}(p^1)$ for randomly generated points over $1000$ trials.} \label{fig9} \vspace{-15pt}
\end{figure*}

Note that if $n = 1$, we know that any randomly generated point must be in $\mathcal{N}_{G}(p^1)$. Therefore, the expected number is $1$.
For $n \geq 2$, we examine it via simulations. We randomly generate $n$ points, run $5000$ trials, and calculate the average number of points in $\mathcal{N}_{G}(p^1)$. The plot of the average number of points in $\mathcal{N}_{G}(p^1)$ as a function of the number of points is given in Fig.~\ref{fig17}. As can be seen from the figure, the average number of points in $\mathcal{N}_{G}(p^1)$ (which is an approximation of the expected number of points in $\mathcal{N}_{G}(p^1)$) always stays below $4.5$. This is also confirmed via simulations for more points (though less trials are run for the sake of time). More specifically, we choose $60, 80, 100, 120, 140, 160$ points, and fix the number of trials to be $1000$. For example, if $|\mathds{P}| = 60$ and we run $1000$ trials, there are $4.34$ points in $\mathcal{N}_{G}(p^1)$ on average, and the histogram of the number of points in $\mathcal{N}_{G}(p^1)$ is plotted in the left upper corner of Fig.~\ref{fig9}. Similarly, if $|\mathds{P}| = 80, 100, 120, 140, 160$ and we still run $1000$ trials, the average number of points in $\mathcal{N}_{G}(p^1)$ is below $4.5$, which is independent of the number of points that we generate. The histograms also have similar shapes for all scenarios. The intuition is that, when an additional point is added, the probability (that the number of points in $\mathcal{N}_{G}(p^1)$ increases) decreases as the number of points increases because it is more likely that the point is dominated by other points or other points dominate this point; therefore, the expected number of points in $\mathcal{N}_{G}(p^1)$ will not blow up and is very likely to stay around a certain value when the number of points is large enough.

\subsection{Methods for Solving Problem~\ref{problem:distributed_calculation}}

As discussed in the previous subsection, the number of points in $\mathcal{N}_{G}(p^1)$ is comparatively much smaller than the number of points in $\mathds{P}$ especially when the total number of points is large. Therefore, it is much more efficient to calculate the Voronoi cell of the point $p^1$ based on the points in $\mathcal{N}_{G}(p^1)$. In this subsection, we briefly discuss methods for solving Problem~\ref{problem:distributed_calculation}.

Given a fixed point $p^1$ and a set of points $\mathds{P}$, one straightforward way to calculate the Voronoi cell of $p^1$ is to consider one point $p \in \mathds{P}$ at a time and calculate the boundary based on the points that have considered so far; we denote this approach as the naive approach. The upper bound based approach is that we first calculate $\mathcal{N}_{G}(p^1)$, and then calculate the boundary of the Voronoi cell only based on the points in $\mathcal{N}_{G}(p^1)$. The difference between these two approaches lies in the fact that the naive approach uses all points to directly calculate the Voronoi cell, and the upper bound based approach only uses points in $\mathcal{N}_{G}(p^1)$ to do the calculation. Since the calculation of the upper bound is much simpler compared with calculating the Voronoi cell of $p^1$ (which involves determining the intersection points of hyperbolas), the additional effort to preselect the set of points that are necessary to compute the Voronoi cell is well worthy.

Intuitively, the ratio of the naive approach's running time to that of the dominance graph based approach will roughly be the ratio of the total number of points to the number of points in $\mathcal{N}_{G}(p^1)$, which is denoted as $\mathcal{R}$, i.e., $\mathcal{R} = \frac{|\mathds{P}|}{|\mathds{N}_G(p^1)|}$.  To calculate the Voronoi cell, we can use any method that is capable of calculating Voronoi cells for additively weighted metrics. For example, the most efficient method is the sweepline algorithm proposed in~\cite{yru_journal:Fortune_1987}; however, the algorithm computes the Voronoi cell for every point instead of just one single point, and the implementation is complicated since queuing mechanism is necessary. In our simulations, we implement a simple (and less efficient) algorithm for calculating the (bounded) Voronoi cell of the point $p^1$. The basic idea is that we consider one point $p$ at a time, and determine which part of the boundary between $p$ and $p^1$ contributes to the (bounded) Voronoi cell of $p^1$.

\begin{figure}[tb] \centering
\subfigure[Point locations.]{\includegraphics[scale=0.28]{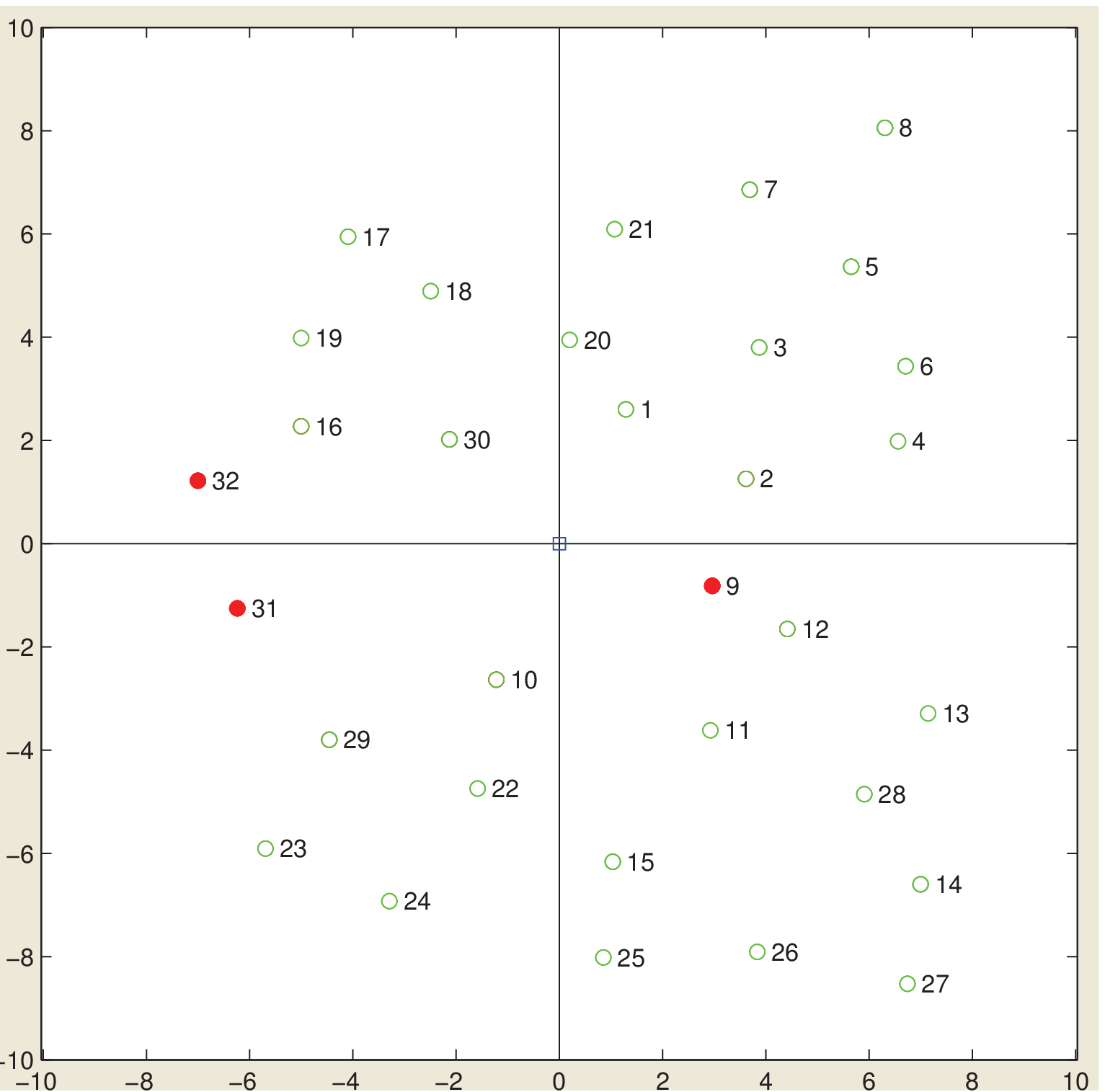}} \hfil
\subfigure[Voronoi cell.]{\includegraphics[scale=0.29]{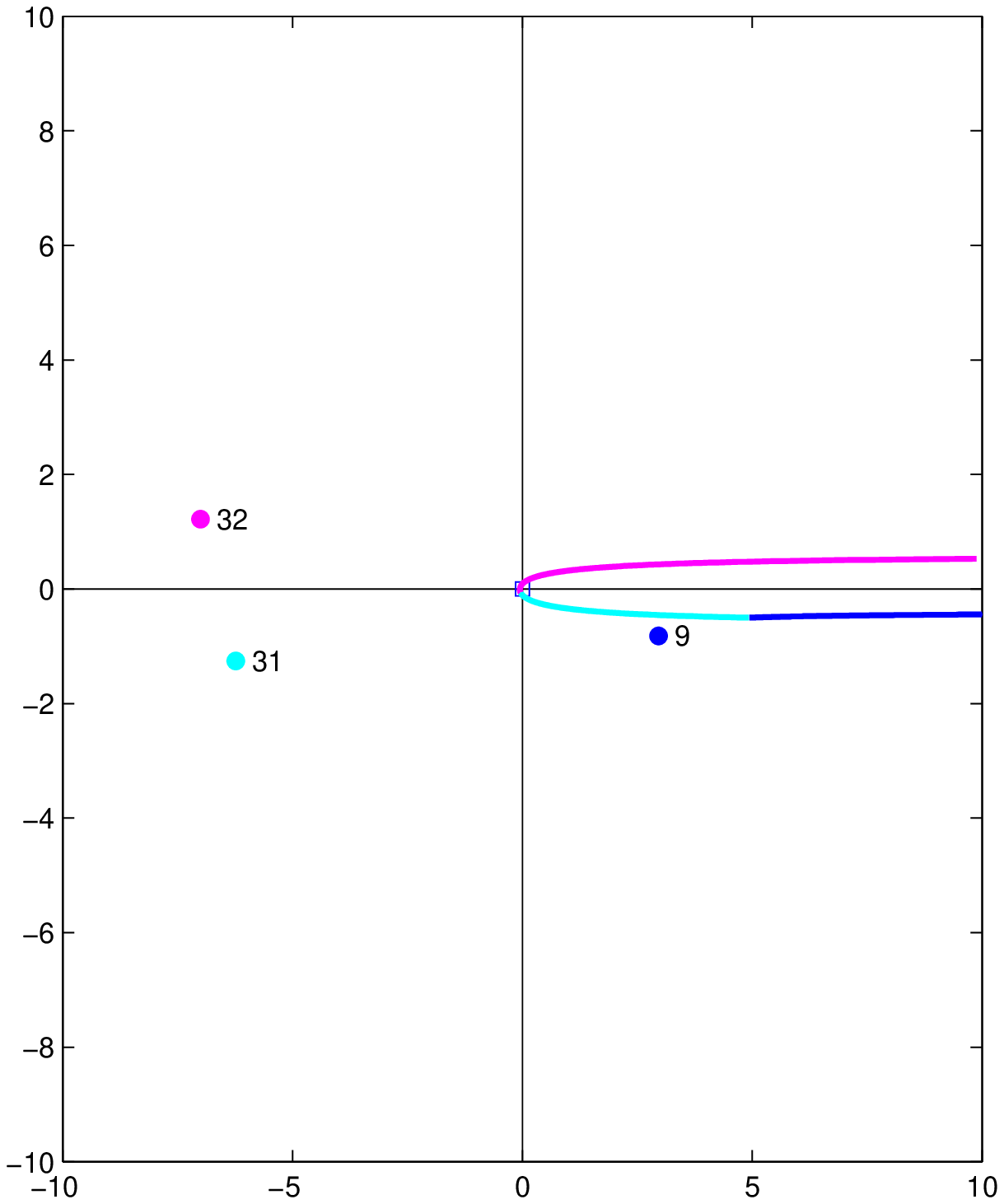}}
\caption{A fixed point $p^1 = (0~0)^T$ and a set of $32$ points in (a), and the (bounded) Voronoi cell of the point $p^1$ in (b).} \label{fig15}
\end{figure}

\begin{figure}[tb] \centering
\psfrag{1}{$\mathcal{R}$}
\includegraphics[scale=0.25]{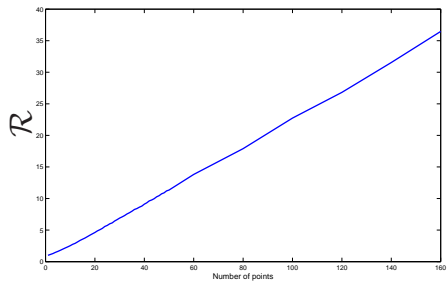}
\caption{Plot of the ratio of the number of points to the average number of points in $\mathcal{N}_{G}(p^1)$.} \label{fig18} \vspace{-15pt}
\end{figure}

We generate $32$ points as shown in Fig.~\ref{fig15}(a), and calculate the Voronoi cell. For the naive approach, $\mathds{P'} = \{1, 2, 3, ..., 32\}$; for the upper bound based approach, $\mathds{P'} = \{9, 31, 32\}$ which is the set of points in $\mathcal{N}_{G}(p^1)$. The Voronoi cells are plotted in Fig.~\ref{fig15}(b). Note that the upper bound based approach runs\footnote{Note that the running time for the upper bound based approach does not include the time necessary to compute $\mathcal{N}_{G}(p^1)$, which is ignorable when the total number of points is large.} $1.55$ seconds while the naive approach runs $18.75$ seconds. The running time ratio is $\frac{18.75}{1.55} \approx 12.1$, while the ratio $\mathcal{R}$ is $\frac{32}{3} \approx 10.7$. Note that the two ratios are close as expected. By examining the Voronoi cells, they are exactly the same. In Fig.~\ref{fig18}, the ratio $\mathcal{R}$ of the number of all points to the average number of points in $\mathcal{N}_{G}(p^1)$ is plotted by combing the data in Fig.~\ref{fig17} and Fig.~\ref{fig9} (note that for the number of points larger than $60$, the plot is based on interpolation using the data in Fig.~\ref{fig9}). This plot shows that on average the ratio $\mathcal{R}$ increases linearly as a function of the number of points. In other words, the upper bound based approach becomes more and more efficient (on average) as the number of points increases. Note that even in the worst case, i.e., the largest number of points in $\mathcal{N}_{G}(p^1)$ given a fixed number of points, the ratio $\mathcal{R}$ is also increasing as can be verified from Fig.~\ref{fig9}.

\section{Conclusions} \label{section8}

In this paper, we study the Voronoi partition by introducing an energy-based metric in constant flow environments. We provide an explicit expression for this energy metric, and use it to derive the equation for the Voronoi boundary between two generators. To facilitate the distributed calculation of Voronoi cells, we propose a disk-based approximation which leads to a lower bound on the set of Voronoi neighbors, and asymptote-based approximations which lead to an upper bound on the set of Voronoi neighbors. When deriving the upper bound, we introduce the dominance relation and provide a complete characterization.  Simulations are run to evaluate the upper bound and its effect on calculating Voronoi cells. The results have potential applications to any other setting based on additively weighted metrics.

There are several future directions. First, we would like to study the Voronoi partition when an upper bound on the traveling time of vehicles is imposed. This is motivated by applications such as search/rescue in which the time taken to reach a point is also very important besides saving energy. Second, we would like to generalize the flow field to more general scenarios, such as piecewise constant flows~\cite{yru_journal:Kwok_2010}, time varying flows~\cite{yru_journal:Bakolas_2010_Voronoi_j}, or even quadratic flows~\cite{yru_journal:Ru_2011_b}. Third, we would also like to incorporate the energy necessary for communication~\cite{yru_journal:Rodoplu_1999} (such as transmitting the location information) and sensing into our energy metric. Last, we would like to study approximation techniques for Voronoi cells based on metrics other than the energy metric. For example, we would like to extend such results to the power metric as studied in~\cite{yru_journal:Pavone_2011}.

\appendix
\noindent \textbf{Proof of Proposition~\ref{prop:metric}}
The minimum energy control problem can be formulated as below:
\begin{align*}
\min ~&\int_{0}^{t_f} U^T U d t \nonumber\\
s.t. ~~& \frac{dx}{dt} = U_x + B~, ~~\frac{dy}{dt} = U_y~,\nonumber\\
& x(0) = x_{p^1}, ~y(0) = y_{p^1}~, x(t_f) = x_{p^2}, ~y(t_f) = y_{p^2}~.
\end{align*}
The objective of the optimization problem is to find a control $U$ which minimizes the total energy. Then the Hamiltonian is $H = U^TU + P^T (U + N)$, where $P= (P_1~P_2)^T$ and $N = (B~0)^T$. Using the minimum principle~\cite{yru_journal:Bryson_1975}, we obtain the
following coupled ordinary differential equations (ODEs) besides Eqs.~\eqref{eq:dynamicsx} and \eqref{eq:dynamicsy}:
\begin{align}
\frac{dP_1}{dt} &= 0~, \label{eq:dP_1}\\
\frac{dP_2}{dt} &= 0~. \label{eq:dP_2}
\end{align}
Since $U$ is chosen to minimize the Hamiltonian, we have $U =
-\frac{1}{2} P$. Plugging $U$ into Eqs.~(\ref{eq:dynamicsx}) and
(\ref{eq:dynamicsy}), we get the following ODEs:
\begin{align}
\frac{dx}{dt} &= -\frac{1}{2}P_1 + B~,\label{eq:dxnew}\\
\frac{dy}{dt} &= -\frac{1}{2} P_2~. \label{eq:dynew}
\end{align}

From Eq.~(\ref{eq:dP_1}), we know $P_1$ is a constant, and let
$P_1(t) = C_1$ for $t \in [0, t_f]$. Similarly, from
Eq.~(\ref{eq:dP_2}), we know $P_2$ is also a constant and let
$P_2(t) = C_2$ for $t \in [0, t_f]$. Then using
Eqs.~(\ref{eq:dxnew}) and (\ref{eq:dynew}), we have
$$x(t) = x_{p^1} + (B -
\frac{C_1}{2})t,~\mathrm{and}~y(t) = y_{p^1} -
\frac{C_2 t}{2}~.$$ Since $x(t_f)$ and $y(t_f)$ are given, we
obtain the following equations
\begin{align}
x_{p^2} &= x_{p^1} + (B - \frac{C_1}{2})t_f~, \label{eq:setofeq_1}\\
y_{p^2} &= y_{p^1} - \frac{C_2 t_f}{2}~. \label{eq:setofeq_2}
\end{align}

Since $t_f$ is free and there is no cost imposed on the final state in the
optimization problem, $H|_{tf} = 0$. Therefore, we have
\begin{equation}
-\frac{C_1^2 + C_2^2}{4} + C_1 B = 0~. \label{eq:setofeq_3}
\end{equation}

By solving Eqs.~(\ref{eq:setofeq_1}), (\ref{eq:setofeq_2}) and
(\ref{eq:setofeq_3}), we obtain the unique solution in Eq.~\eqref{eq:tf}. Then the optimal control is $U(t) = -\frac{1}{2}
\begin{bmatrix} C_1\\C_2\end{bmatrix}$ for $t \in [0, t_f]$, and the minimum energy is given in Eq.~\eqref{eq:min_energy}.\\

\noindent \textbf{Proof of Theorem~\ref{theorem:Voronoi_neighbor_lower}}
We prove it by contradiction. Suppose $p^k \in \mathcal{N}_D(p^1)$, but $p^k \notin \mathcal{N}_V(p^1)$. Then the point $p^*(p^1, p^k)$ must lie outside of $V(p^1)$. In other words, there must exist an edge of the Voronoi cell $V(p^1)$ that intersects with the line segment $p^1 p^*(p^1, p^k)$ at a point $p^l$ (which is different from $p^*(p^1, p^k)$) satisfying $d_{p^1 p^l} < d_{p^1 p^*(p^1, p^k)}$. Suppose the edge is due to another point $p^m$. Since $p^*(p^1, p^m)$ is the unique point on the boundary between $p^1$ and $p^m$ that is closest to $p^1$, we have $d_{p^1, p^*(p^1, p^m)} \leq d_{p^1 p^l}$. Thus, $d_{p^1, p^*(p^1, p^m)} \leq d_{p^1 p^l} < d_{p^1 p^*(p^1, p^k)}$. That is to say, following the disk-based approximation, $d_{p^1, p^*(p^1, p^m)}$ (i.e., the radius due to the point $p^m$) is strictly smaller than $d_{p^1, p^*(p^1, p^k)}$ (i.e., the radius due to the point $p^k$). It contradicts with the assumption that $p^k \in \mathcal{N}_D(p^1)$. Therefore, $p^k$ must be a Voronoi neighbor.\\

\noindent \textbf{Proof of Proposition~\ref{prop:neighbor_lower}}
Since $p^*(p^1, p^k)$ lies on the boundary between $p^1$ and $p^k$, we have $J(p^1, p^*(p^1, p^k)) = J(p^k, p^*(p^1, p^k))$. Because for any $l \in \{2, 3, ..., n\} \setminus \{k\}$ we have $J(p^k, p^*(p^1, p^k)) < J(p^l, p^*(p^1, p^k))$, we obtain $J(p^1, p^*(p^1, p^k)) \leq J(p^m, p^*(p^1, p^k))$ for $m \in \{2, 3, ..., n\}$. Based on Definition~\ref{def:Voronoi}, we know that $p^*(p^1, p^k)$ lies in $V(p^1)$. Because for any $l \in \{2, 3, ..., n\} \setminus \{k\}$ we have $J(p^1, p^*(p^1, p^k)) < J(p^l, p^*(p^1, p^k))$, $p^*(p^1, p^k)$ must only belong to the edge of $V(p^1)$ between $p^1$ and $p^k$, and the edge of $V(p^1)$ between $p^1$ and $p^k$ is not a single point due to the strict inequality. Therefore, $p^k$ must be a Voronoi neighbor of $p^1$.\\

\noindent \textbf{Proof of Proposition~\ref{prop:asymptotes}}
For $l_1$, we are interested in $y' = \frac{b}{a}x'$ for $y' \geq 0$. Plugging in the expressions for $x'$ and $y'$ in Eqs.~\eqref{eq:x_new} and~\eqref{eq:y_new}, we obtain
\begin{align*}
&-(x - \frac{x_{p^1} + x_{p^2}}{2}) \sin \alpha + (y -
\frac{y_{p^1} + y_{p^2}}{2}) \cos \alpha =\\
 &((x - \frac{x_{p^1} + x_{p^2}}{2}) \cos \alpha + (y - \frac{y_{p^1} + y_{p^2}}{2}) \sin \alpha) \tan \alpha~,
\end{align*}
which simplifies to be $y - \frac{y_{p^1} + y_{p^2}}{2} = (x - \frac{x_{p^1} + x_{p^2}}{2}) \tan2 \alpha$. $y' \geq 0$ is equivalent to $-(x - \frac{x_{p^1} + x_{p^2}}{2}) \sin \alpha + (y -
\frac{y_{p^1} + y_{p^2}}{2}) \cos \alpha \geq 0$. Using $x - \frac{x_{p^1} + x_{p^2}}{2} = \frac{y - \frac{y_{p^1} + y_{p^2}}{2}}{\tan2 \alpha}$, we obtain $(y - \frac{y_{p^1} + y_{p^2}}{2}) \frac{\sin \alpha}{\sin 2 \alpha} > 0$. Since $\alpha \in (0, \frac{\pi}{2})$, $y' \geq 0$ reduces to $y \geq \frac{y_{p^1} + y_{p^2}}{2}$.

For $l_2$, we are interested in $y' = -\frac{b}{a}x'$ for $y' \leq 0$. Plugging in the expressions for $x'$ and $y'$ in Eqs.~\eqref{eq:x_new} and~\eqref{eq:y_new}, we obtain
\begin{align*}
&-(x - \frac{x_{p^1} + x_{p^2}}{2}) \sin \alpha + (y -
\frac{y_{p^1} + y_{p^2}}{2}) \cos \alpha =\\
 &((x - \frac{x_{p^1} + x_{p^2}}{2}) \cos \alpha + (y - \frac{y_{p^1} + y_{p^2}}{2}) \sin \alpha) (-\tan \alpha)~,
\end{align*}
which can be rewritten as $(y - \frac{y_{p^1} + y_{p^2}}{2})\frac{1}{\cos \alpha} = 0$. Since $\alpha \in (0, \frac{\pi}{2})$, we have $y = \frac{y_{p^1} + y_{p^2}}{2}$. $y' \leq 0$ is equivalent to $-(x - \frac{x_{p^1} + x_{p^2}}{2}) \sin \alpha \leq 0$. Therefore, $x \geq \frac{x_{p^1} + x_{p^2}}{2}$.\\

\noindent \textbf{Proof of Theorem~\ref{case_1}}

If $p^3$ is the same as $p^2$, $p^2 \succ p^3$ and the conditions hold trivially. Therefore, in the following proof, we assume that $p^2$ and $p^3$ are different.

\noindent \textbf{(If part)} Depending on the $x$ coordinates of $p^1$ and $p^2$, there are three cases.

\begin{figure}[tb] \centering

\psfrag{0}{$0$} \psfrag{x}{$x$} \psfrag{y}{$y$} \psfrag{1}{$p^1$} \psfrag{2}{$q_1$} \psfrag{3}{$\alpha$} \psfrag{4}{$p^2$}
\psfrag{5}{$q_2$} \psfrag{6}{$q_3$} \psfrag{7}{$p^3$} \psfrag{8}{$q_4$} \psfrag{9}{$q_5$} \psfrag{10}{$q_6$} \psfrag{11}{$\beta$}
\includegraphics[scale=0.3]{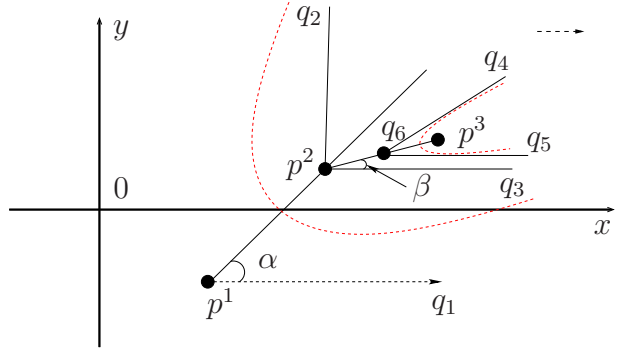}
\caption{Points $\{p^1, p^2\}$ satisfying $x_{p^1} < x_{p^2}$ and $y_{p^1} < y_{p^2}$. The red dashed lines are the boundaries.} \label{fig5} \vspace{-15pt}
\end{figure}

\textbf{Case I}: $x_{p^2} > x_{p^1}$. Since $y_{p^3} \geq y_{p^2} > y_{p^1}$, we have $(x_{p^3} - x_{p^1}) \times (y_{p^2} - y_{p^1}) \geq (y_{p^3} - y_{p^1}) \times (x_{p^2} - x_{p^1}) > 0$, which implies that $x_{p^3} > x_{p^1}$. Therefore, Eq.~\eqref{eq:inequality} can be rewritten as
\begin{equation}
\frac{y_{p^3} - y_{p^1}}{x_{p^3} - x_{p^1}} \leq \frac{y_{p^2} - y_{p^1}}{x_{p^2} - x_{p^1}}~. \label{eq:angle_relation}
\end{equation}
Now there are two cases depending on the $y$ coordinate of $p^2$ and $p^3$:
\begin{itemize}
\item If $y_{p^3} = y_{p^2}$, Eq.~\eqref{eq:angle_relation} implies that $x_{p^3} > x_{p^2}$ since $p^2, p^3$ are different. For any point $p$ in $V(p^1~|~p^2)$, we have $J(p^1, p) \leq J(p^2, p)$, i.e., $2B(d_{p^1p} + x_{p^1} - x_p) \leq 2B(d_{p^2p} + x_{p^2} - x_p)$. Since $d_{p^2p} \leq d_{p^3p} + d_{p^2p^3} = d_{p^3p} + x_{p^3} - x_{p^2}$ (due to the triangular inequality), i.e., $d_{p^2p} + x_{p^2} \leq d_{p^3p} + x_{p^3}$, we have $2B(d_{p^2p} + x_{p^2} - x_p) \leq 2B(d_{p^3p} + x_{p^3} - x_p)$. Therefore, $2B(d_{p^1p} + x_{p^1} - x_p) \leq 2B(d_{p^3p} + x_{p^3} - x_p)$, i.e., $J(p^1, p) \leq J(p^3 p)$. Therefore, we have $V(p^1~|~p^2) \subseteq V(p^1~|~p^2, p^3)$. As $V(p^1~|~p^2, p^3) \subseteq V(p^1~|~p^2)$, we have $V(p^1~|~p^2) = V(p^1~|~p^2, p^3)$;

\item If $y_{p^3} > y_{p^2}$, Eq.~\eqref{eq:angle_relation} implies that $x_{p^3} > x_{p^2}$. For any point $p$ in $V(p^1~|~p^2)$, we have $p \in \subscr{D}{upper}''(p^1~|~p^2)$, where $\subscr{D}{upper}''(p^1~|~p^2)$ is obtained via shifting the vertex of $\subscr{D}{upper}(p^1~|~p^2)$ to the point $p^2$. In other words, the two lines that consist of the boundary of $\subscr{D}{upper}''(p^1~|~p^2)$ are $y = y_{p^2}$ with $x \geq x_{p^2}$ (namely, the half line $p^2 q_3$ in Fig.~\ref{fig5}), and $y - y_{p^2} = (x - x_{p^2}) \tan 2 \alpha$ with $y \geq y_{p^2}$ and $\alpha$ being the angle of $\angle q_1p^1p^2$ (namely, the half line $p^2 q_2$ in Fig.~\ref{fig5}). Now we look at $\subscr{D}{lower}(p^2~|~p^3)$, the lower approximation for $p^2$ given $p^3$. Recall that the two lines that consist of the boundary of $\subscr{D}{lower}(p^2~|~p^3)$ are $y = \frac{y_{p^2} + y_{p^3}}{2}$ with $x \geq \frac{x_{p^2} + x_{p^3}}{2}$ (namely, the half line $q_6 q_5$ in Fig.~\ref{fig5}, where $q_6$ is the middle point of $p^2, p^3$), and $y - \frac{y_{p^2} + y_{p^3}}{2} = (x - \frac{x_{p^2} + x_{p^3}}{2}) \tan 2 \beta$ with $y \geq \frac{y_{p^2} + y_{p^3}}{2}$ and $\beta$ being the angle of $\angle q_3p^2p^3$ (namely, the half line $q_6 q_4$ in Fig.~\ref{fig5}). Since $\frac{y_{p^2} + y_{p^3}}{2} > y_{p^2}$, $\frac{x_{p^2} + x_{p^3}}{2} > x_{p^2}$ and $\beta \leq \alpha$ (due to Eq.~\eqref{eq:angle_relation}), we have $2 \beta \leq 2 \alpha$, and $\subscr{D}{upper}''(p^1~|~p^2) \subseteq \subscr{D}{lower}(p^2~|~p^3)$. In other words, for any $p \in V(p^1~|~p^2)$, $p \in \subscr{D}{lower}(p^2~|~p^3)$. Therefore, $2B(d_{p^1p} + x_{p^1} - x_p) \leq 2B(d_{p^2p} + x_{p^2} - x_p) \leq 2B(d_{p^3p} + x_{p^3} - x_p)$. Thus, $V(p^1~|~p^2) \subseteq V(p^1~|~p^2, p^3)$. Since $V(p^1~|~p^2, p^3) \subseteq V(p^1~|~p^2)$, we have $V(p^1~|~p^2) = V(p^1~|~p^2, p^3)$.
\end{itemize}

\textbf{Case II}: $x_{p^2} = x_{p^1}$. Since $y_{p^2} > y_{p^1}$, Eq.~\eqref{eq:inequality} implies that $x_{p^3} \geq x_{p^1} = x_{p^2}$. There are two cases depending on the $x$ coordinates of $p^2, p^3$:
\begin{itemize}
\item If $x_{p^3} = x_{p^2}$, then we have $x_{p^3} = x_{p^1}$ and $y_{p^3} > y_{p^2}$ since $p^2, p^3$ are different. In this case, the boundary between $p^1$ and $p^2$ is $y_{12} = \frac{y_{p^1} + y_{p^2}}{2}$ as shown in Case II of Section~\ref{section3_special}, and the boundary between $p^1$ and $p^3$ is $y_{13} = \frac{y_{p^1} + y_{p^3}}{2} > y_{12}$. Therefore, we have $V(p^1~|~p^2) = V(p^1~|~p^2, p^3)$;

\item If $x_{p^3} > x_{p^2}$, then we have $x_{p^3} > x_{p^1}$ and $y_{p^3} \geq y_{p^2}$. In this case, the boundary between $p^1$ and $p^2$ is $y_{12} = \frac{y_{p^1} + y_{p^2}}{2}$, and the boundary between $p^1$ and $p^3$ is a hyperbola which is above the lower asymptote $y_{13} = \frac{y_{p^1} + y_{p^3}}{2} \geq y_{12}$. However, the hyperbola will not touch the line $y_{12}$. Therefore, we have $V(p^1~|~p^2) = V(p^1~|~p^2, p^3)$.
\end{itemize}

\textbf{Case III}: $x_{p^2} < x_{p^1}$. There are three cases depending on the $x$ coordinates of $p^2, p^3$:
\begin{itemize}
\item If $x_{p^3} > x_{p^1}$, then Eq.~\eqref{eq:inequality} holds because $y_{p^3} \geq y_{p^2} > y_{p^1}$. The boundary between $p^1$ and $p^2$ is a hyperbola which is below the upper asymptote $y_{12} = \frac{y_{p^1} + y_{p^2}}{2}$, while the boundary between $p^1$ and $p^3$ is a hyperbola which is above the lower asymptote $y_{13} = \frac{y_{p^1} + y_{p^3}}{2} \geq y_{12}$. Since the hyperbolas will not touch the asymptotes, we have $V(p^1~|~p^2) = V(p^1~|~p^2, p^3)$;

\item If $x_{p^3} = x_{p^1}$, the boundary between $p^1$ and $p^2$ is a hyperbola which is below the upper asymptote $y_{12} = \frac{y_{p^1} + y_{p^2}}{2}$, while the boundary between $p^1$ and $p^3$ is the line $y_{13} = \frac{y_{p^1} + y_{p^3}}{2} \geq y_{12}$. Since the hyperbola will not touch the asymptote $y_{12}$ due to $y_{p^1} < y_{p^2}$, we have $V(p^1~|~p^2) = V(p^1~|~p^2, p^3)$;

\item If $x_{p^3} < x_{p^1}$, it can be shown that $V(p^1~|~p^2) = V(p^1~|~p^2, p^3)$ by an argument similar to the one used in Case I.
\end{itemize}

\noindent \textbf{(Only if part)} We prove it via contradiction. Suppose the conditions that $y_{p^3} \geq y_{p^2}$ and $(y_{p^3} - y_{p^1}) \times (x_{p^2} - x_{p^1}) \leq (x_{p^3} - x_{p^1}) \times (y_{p^2} - y_{p^1})$ do not hold. Then there are three possibilities: 1) $y_{p^3} < y_{p^1}$, 2) $y_{p^1} \leq y_{p^3} < y_{p^2}$, and 3) $y_{p^3} \geq y_{p^2}$ but $(y_{p^3} - y_{p^1}) \times (x_{p^2} - x_{p^1}) > (x_{p^3} - x_{p^1}) \times (y_{p^2} - y_{p^1})$.

\textbf{Case 1)}: $y_{p^3} < y_{p^1}$. In this case, $y_{p^3} < y_{p^1} < y_{p^2}$. The boundary between $p^1$ and $p^3$ is a hyperbola which is below the upper asymptote $y_{13} = \frac{y_{p^1} + y_{p^3}}{2} < y_{p^1}$. Since $y_{p^2} > y_{p^1}$, the boundary between $p^1$ and $p^2$ is a hyperbola which is above the lower asymptote $y_{12} = \frac{y_{p^1} + y_{p^2}}{2} > y_{p^1}$, and lies strictly above the boundary between $p^1$ and $p^3$. Therefore, we must have $V(p^1~|~p^2, p^3) \subset V(p^1~|~p^2)$. A contradiction to $p^2 \succ p^3$.

\textbf{Case 2)}: $y_{p^1} \leq y_{p^3} < y_{p^2}$. Then there are two cases depending on the $x$ coordinates of $p^1, p^2$:
\begin{itemize}
\item $x_{p^2} \geq x_{p^1}$. In this case, the boundary between $p^1$ and $p^2$ lies on/above the asymptote $y_{12} = \frac{y_{p^1} + y_{p^2}}{2}$. If $x_{p^3} \geq x_{p^1}$ (or $x_{p^3} < x_{p^1}$), the boundary between $p^1$ and $p^3$ lies on/above (or below) the asymptote $y_{13} = \frac{y_{p^1} + y_{p^3}}{2} < y_{12}$ and can be arbitrarily close to the asymptote. Therefore, we must have $V(p^1~|~p^2, p^3) \subset V(p^1~|~p^2)$. A contradiction.

\item $x_{p^2} < x_{p^1}$. In this case, the boundary between $p^1$ and $p^2$ lies below the asymptote $y_{12} = \frac{y_{p^1} + y_{p^2}}{2}$ and can be arbitrarily close to the asymptote. If $x_{p^3} \geq x_{p^1}$ (or $x_{p^3} < x_{p^1}$), the boundary between $p^1$ and $p^3$ lies on/above (or below) the asymptote $y_{13} = \frac{y_{p^1} + y_{p^3}}{2} < y_{12}$ and can be arbitrarily close to the asymptote. Therefore, we must have $V(p^1~|~p^2, p^3) \subset V(p^1~|~p^2)$. A contradiction.
\end{itemize}

\textbf{Case 3)}: $y_{p^3} \geq y_{p^2}$ but $(y_{p^3} - y_{p^1}) \times (x_{p^2} - x_{p^1}) > (x_{p^3} - x_{p^1}) \times (y_{p^2} - y_{p^1})$. There are also two cases depending on the $x$ coordinates of $p^1, p^2$:
\begin{itemize}
\item $x_{p^2} \geq x_{p^1}$. The boundary between $p^1$ and $p^2$ lies on/above the asymptote $y_{12} = \frac{y_{p^1} + y_{p^2}}{2}$, and can be arbitrarily close to the asymptote. If $x_{p^3} < x_{p^1}$, the boundary between $p^1$ and $p^3$ lies below the asymptote $y_{13} = \frac{y_{p^1} + y_{p^3}}{2} \geq y_{12}$ and can be arbitrarily close to the asymptote. Therefore, we must have $V(p^1~|~p^2, p^3) \subset V(p^1~|~p^2)$. A contradiction. If $x_{p^3} \geq x_{p^1}$, the upper asymptote of the boundary between $p^1$ and $p^3$ has the slope $\tan 2 \beta$ with $\beta = \arctan \frac{y_{p^3} - y_{p^1}}{x_{p^3} - x_{p^1}} \geq 0$. The upper asymptote of the boundary between $p^1$ and $p^2$ has the slope $\tan 2 \alpha$ with $\alpha = \arctan \frac{y_{p^2} - y_{p^1}}{x_{p^2} - x_{p^1}} \geq 0$. If $(y_{p^3} - y_{p^1}) \times (x_{p^2} - x_{p^1}) > (x_{p^3} - x_{p^1}) \times (y_{p^2} - y_{p^1})$, we have $\beta > \alpha$, which implies that $2 \beta > 2 \alpha$. Therefore, part of the boundary (corresponding to the upper asymptote) between $p^1$ and $p^3$ must strictly refine $V(p^1|p^2)$, i.e., $V(p^1~|~p^2, p^3) \subset V(p^1~|~p^2)$. A contradiction.

\item $x_{p^2} < x_{p^1}$. Since $y_{p^3} \geq y_{p^2} > y_{p^1}$, $(y_{p^3} - y_{p^1}) \times (x_{p^2} - x_{p^1}) > (x_{p^3} - x_{p^1}) \times (y_{p^2} - y_{p^1})$ implies that $x_{p^3} < x_{p^1}$. The lower asymptote of the boundary between $p^1$ and $p^3$ has the slope $\tan 2 \beta$ with $\beta = \arctan \frac{y_{p^3} - y_{p^1}}{x_{p^3} - x_{p^1}} \leq 0$. The lower asymptote of the boundary between $p^1$ and $p^2$ has the slope $\tan 2 \alpha$ with $\alpha = \arctan \frac{y_{p^2} - y_{p^1}}{x_{p^2} - x_{p^1}} \leq 0$. If $(y_{p^3} - y_{p^1}) \times (x_{p^2} - x_{p^1}) > (x_{p^3} - x_{p^1}) \times (y_{p^2} - y_{p^1})$, we have $\beta > \alpha$, i.e., $0 \leq -\beta < -\alpha$, which implies that $0 \leq -2 \beta <- 2 \alpha$. Therefore, part of the boundary (corresponding to the lower asymptote) between $p^1$ and $p^3$ must strictly refine $V(p^1|p^2)$, i.e., $V(p^1~|~p^2, p^3) \subset V(p^1~|~p^2)$. A contradiction.
\end{itemize}

\vspace{15pt} \noindent \textbf{Proof of Proposition~\ref{prop:antisymmetry}}
If $p^2$ is the same as $p^3$, the result holds trivially. In the following proof, we consider the case in which $p^2$ and $p^3$ are different.

If $y_{p^1} < y_{p^2}$ (namely, Scenario A in Section~\ref{section5_scenario}), $p^2 \succ p^3$ implies that $y_{p^3} \geq y_{p^2}$, and $(y_{p^3} - y_{p^1}) \times (x_{p^2} - x_{p^1}) \leq (x_{p^3} - x_{p^1}) \times (y_{p^2} - y_{p^1})$. Since $y_{p^3} \geq y_{p^2} > y_{p^1}$, $p^3 \succ p^2$ implies that $y_{p^2} \geq y_{p^3}$, and $(y_{p^2} - y_{p^1}) \times (x_{p^3} - x_{p^1}) \leq (x_{p^2} - x_{p^1}) \times (y_{p^3} - y_{p^1})$. From $y_{p^3} \geq y_{p^2}$ and $y_{p^2} \geq y_{p^3}$, we have $y_{p^2} = y_{p^3} > y_{p^1}$. Therefore, $(y_{p^3} - y_{p^1}) \times (x_{p^2} - x_{p^1}) \leq (x_{p^3} - x_{p^1}) \times (y_{p^2} - y_{p^1})$ implies $x_{p^2} \leq x_{p^3}$, and $(y_{p^2} - y_{p^1}) \times (x_{p^3} - x_{p^1}) \leq (x_{p^2} - x_{p^1}) \times (y_{p^3} - y_{p^1})$ implies $x_{p^3} \leq x_{p^2}$. Thus, $x_{p^2} = x_{p^3}$. In summary, $p^2 = p^3$. Similarly, we can show that $p^2 = p^3$ for the case $y_{p^1} > y_{p^2}$ (namely, Scenario B in Section~\ref{section5_scenario}).

If $y_{p^1} = y_{p^2}$ and $x_{p^1} < x_{p^2}$ (namely, Scenario D in Section~\ref{section5_scenario}), $p^2 \succ p^3$ implies that $y_{p^3} = y_{p^2}$ and $x_{p^3} \geq x_{p^2}$. Since $y_{p^3} = y_{p^1}$ and $x_{p^3} \geq x_{p^2} > x_{p^1}$, $p^3 \succ p^2$ implies that $y_{p^2} = y_{p^3}$ and $x_{p^2} \geq x_{p^3}$. Therefore, $p^2 = p^3$.\\

\noindent \textbf{Proof of Theorem~\ref{prop:Voronoi_neighbor}}
Assumption~\ref{assumption} guarantees that $V(p^1)$ is nonempty. Let $p^i \in \mathcal{N}_V(p^1)$ for $i \in \{2, 3, ..., n\}$, i.e., $p^i$ is a Voronoi neighbor of $p^1$. Then the intersection of the Voronoi cells of $p^1$ and $p^i$ is a curve with the set of points being $V(p^1) \cap V(p^i)$. For any $p$ which lies on the curve but is not an end point, we have $J(p^1, p) = J(p^i, p)$, and $J(p^1, p) < J(p^j, p)$ for any $j \in \{2, 3, ..., n\} \setminus \{i\}$.

Suppose $p^i \notin \mathcal{N}_G(p^1)$. Then there exists some $p^k$ for $k \in \{2, 3, ..., n\} \setminus \{i\}$ such that $p^i$ is dominated by $p^k$. Since $p^k \succ p^i$, $V(p^1|p^k) = V(p^1|p^k, p^i)$. Since for the $p$ chosen previously $p \in V(p^1) \subseteq V(p^1 | p^k)$, we have $J(p^1, p) \leq J(p^k, p)$. Since $V(p^1|p^k) = V(p^1|p^k, p^i)$, $J(p^1, p) \leq J(p^k, p)$ implies that $J(p^1, p) \leq J(p^i, p)$. Therefore, we must have $J(p^1, p) \leq J(p^k, p) \leq J(p^i, p)$ (this can be proved via contradiction). Because $J(p^1, p) = J(p^i, p)$, we have $J(p^1, p) = J(p^k, p)$. This contradicts with $J(p^1, p) < J(p^k, p)$ due to $J(p^1, p) < J(p^j, p)$ for any $j \in \{2, 3, ..., n\} \setminus \{i\}$. Therefore, $\mathcal{N}_V(p^1) \subseteq \mathcal{N}_G(p^1)$.\\

\noindent \textbf{Proof of Proposition~\ref{prop:transitivity}}
If any pair of points in $p^2, p^3, p^4$ is the same, or all three points are the same, the result holds trivially. In the following proof, we consider the case in which all three points are distinct.

Depending on the relative position of $p^1$ and $p^2$, there are four scenarios as discussed in Section~\ref{section5_scenario}:

\textbf{Scenario A} In this case, $y_{p^1} < y_{p^2}$. Since $p^2 \succ p^3$ and $p^3 \succ p^4$, we have $y_{p^1} < y_{p^2} \leq y_{p^3} \leq y_{p^4}$. To show $p^2 \succ p^4$, we only need to prove that
\begin{equation}
(y_{p^4} - y_{p^1}) \times (x_{p^2} - x_{p^1}) \leq (x_{p^4} - x_{p^1}) \times (y_{p^2} - y_{p^1})~. \label{eq:p4p2}
\end{equation}

If $x_{p^1} < x_{p^2}$, $p^2 \succ p^3$ implies that $x_{p^1} < x_{p^3}$ as argued in the proof of Theorem~\ref{case_1}, and $\frac{y_{p^3} - y_{p^1}}{x_{p^3} - x_{p^1}} \leq \frac{y_{p^2} - y_{p^1}}{x_{p^2} - x_{p^1}}$; $p^3 \succ p^4$ implies that $x_{p^1} < x_{p^4}$, and $\frac{y_{p^4} - y_{p^1}}{x_{p^4} - x_{p^1}} \leq \frac{y_{p^3} - y_{p^1}}{x_{p^3} - x_{p^1}}$. Therefore, we have $\frac{y_{p^4} - y_{p^1}}{x_{p^4} - x_{p^1}} \leq \frac{y_{p^2} - y_{p^1}}{x_{p^2} - x_{p^1}}$, which implies Eq.~\eqref{eq:p4p2}.

If $x_{p^1} = x_{p^2}$, $x_{p^3} \geq x_{p^2} = x_{p^1}$. If $x_{p^3} = x_{p^1}$, we must have $x_{p^4} \geq x_{p^1}$ since $p^3 \succ p^4$. Then the left hand side of Eq.~\eqref{eq:p4p2} is $0$ while the right hand side of Eq.~\eqref{eq:p4p2} is nonnegative because $x_{p^4} \geq x_{p^1}$ and $y_{p^2} > y_{p^1}$. Therefore, Eq.~\eqref{eq:p4p2} holds. If $x_{p^3} > x_{p^1}$, then we have $x_{p^4} > x_{p^1}$. Then Eq.~\eqref{eq:p4p2} holds too.

If $x_{p^1} > x_{p^2}$, there are three cases depending on the $x$ coordinates of $p^1, p^3$. i) If $x_{p^3} > x_{p^1}$, then we have $x_{p^4} > x_{p^1}$. Then the left hand side of Eq.~\eqref{eq:p4p2} is negative while the right hand side of Eq.~\eqref{eq:p4p2} is positive because $x_{p^4} > x_{p^1}$ and $y_{p^2} > y_{p^1}$. Therefore, Eq.~\eqref{eq:p4p2} holds. ii) If $x_{p^3} = x_{p^1}$, we must have $x_{p^4} \geq x_{p^1}$. Then Eq.~\eqref{eq:p4p2} holds. iii) If $x_{p^3} < x_{p^1}$ and $x_{p^4} \geq x_{p^1}$, then Eq.~\eqref{eq:p4p2} holds. If $x_{p^3} < x_{p^1}$ and $x_{p^4} < x_{p^1}$, we have $\frac{y_{p^3} - y_{p^1}}{x_{p^3} - x_{p^1}} \leq \frac{y_{p^2} - y_{p^1}}{x_{p^2} - x_{p^1}}$ and $\frac{y_{p^4} - y_{p^1}}{x_{p^4} - x_{p^1}} \leq \frac{y_{p^3} - y_{p^1}}{x_{p^3} - x_{p^1}}$. Therefore, we have $\frac{y_{p^4} - y_{p^1}}{x_{p^4} - x_{p^1}} \leq \frac{y_{p^2} - y_{p^1}}{x_{p^2} - x_{p^1}}$, which implies Eq.~\eqref{eq:p4p2} holds.

\textbf{Scenario B} In this case, $y_{p^1} > y_{p^2}$. It can be proved in a way similar to the one used for Scenario A.

\textbf{Scenario C} In this case, $y_{p^1} = y_{p^2}$ and $x_{p^1} > x_{p^2}$. If $y_{p^2} \neq y_{p^3}$, then $y_{p^3} \neq y_{p^1}$ and the relative position of $p^1$ and $p^3$ belongs to either Scenario A or B. However, in both cases, $p^3 \succ p^4$ implies that $y_{p^4} \neq y_{p^1}$, i.e., $y_{p^4} \neq y_{p^2}$. Then $p^2 \succ p^4$. If $y_{p^2} = y_{p^3}$ and $x_{p^3} < x_{p^1}$, there are two cases:
\begin{itemize}
\item $y_{p^3} \neq y_{p^4}$. Then $y_{p^4} \neq y_{p^2}$. Therefore, $p^2 \succ p^4$.

\item $y_{p^3} = y_{p^4}$ and $x_{p^4} < x_{p^1}$. Then we have $y_{p^2} = y_{p^4}$. Therefore, $p^2 \succ p^4$.
\end{itemize}

\textbf{Scenario D} In this case, $y_{p^1} = y_{p^2}$ and $x_{p^1} < x_{p^2}$.  We have $y_{p^3} = y_{p^2}$ and $x_{p^3} \geq x_{p^2}$ due to $p^2 \succ p^3$, and $y_{p^4} = y_{p^3}$ and $x_{p^4} \geq x_{p^3}$ due to $p^3 \succ p^4$. Therefore, $p^2 \succ p^4$. \\

\noindent \textbf{Proof of Proposition~\ref{prop:acyclic}}
We show the result via contradiction. Suppose there is a cycle $p^i \succ p^{i+1} \succ p^{i+2} ...\succ p^{i+j} \succ p^i$ in the dominance graph where all points are in $\mathds{P}$. Then we have $p^i \succ p^{i+j}$ by repeatedly applying the transitivity property (i.e., Proposition~\ref{prop:transitivity}). Since we also have $p^{i+j} \succ p^i$, $p^i = p^{i+j}$ due to the antisymmetry property (i.e., Proposition~\ref{prop:antisymmetry}). Now the cycle becomes $p^i \succ p^{i+1} \succ p^{i+2} ...\succ p^{i+j-1} \succ p^i$. By repeatedly applying the above analysis, we can show that $p^{i+j-1} = p^{i+j-2} = ... = p^{i+1} = p^i$. Since the points are assumed to be different from each other, a contradiction.


%


%
%


\ifCLASSOPTIONcaptionsoff
  \newpage
\fi




\bibliographystyle{IEEEtran}
\bibliography{IEEEabrv,yru_journal}
\end{document}